\documentclass[letterpaper,11pt]{article}

\newcommand{\an}[2]{{\scriptsize\color{red}}{\color{blue}\  #2}}

\usepackage{endnotes}
\let\footnote=\endnote

\def\bS{\mathbf{S}}
\usepackage{epsfig}
\usepackage{calc}
\usepackage{amstext}
\usepackage{amsmath}
\usepackage{multicol}
\usepackage{amsfonts}
\usepackage{amssymb}
\usepackage{mathrsfs}
\usepackage{bm}
\usepackage{xcolor}
\usepackage{url}
\usepackage{hyperref}
\def\Opt{{\mathrm{Opt}}}
\newcommand{\AR}[2]{\left[\begin{array}{#1}#2\end{array}\right]}

\newcommand{\Tr}{\mathrm{Tr}}
\newcommand{\diag}{\mathrm{diag}}

\newcommand{\rank}{\mathrm{rank}}
\newcommand{\T}{\mathrm{T}}

\newcommand{\cl}{\mathrm{cl}}

\newtheorem{theorem}{Theorem}
\newtheorem{lemma}[theorem]{Lemma}
\newtheorem{corollary}[theorem]{Corollary}
\newtheorem{proposition}[theorem]{Proposition}

\usepackage{epsfig}
\usepackage{amstext}
\usepackage{amsfonts}
\usepackage{amssymb}
\usepackage{mathrsfs}
\usepackage{bm}
\usepackage{cleveref}

\def\bfZ{\mathbf{Z}}
\def\bfE{\mathbf{E}}
\def\Fro{{\mathrm{\tiny Fro}}}
\def\bfzeta{\mathbf{\zeta}}
\def\bfeps{\mathbf{\epsilon}}
\def\bfw{\mathbf{w}}
\def\bfW{\mathbf{W}}
\def\bbE{\mathbb{E}}
\def\mZ{\mathscr{D}[\rho]}

\def\bbR{\mathbb{R}}
\def\bbS{\mathbb{S}}
\def\cT{\mathcal{T}}
\def\Diag{{\mathrm{Diag}}}
\def\cX{{\mathcal{X}}}
\def\cZ{{\mathcal{Z}}}
\def\cA{{\mathcal{A}}}
\def\cF{{\mathcal{F}}}
\def\cS{{\mathcal{S}}}
\def\bK{{\mathbf{K}}}
\def\inter{\mathrm{int}}
\title{
Convex optimization for	finite horizon robust covariance control of linear stochastic systems}
\author{
	Georgios Kotsalis, Guanghui Lan, Arkadi Nemirovski
	\thanks{H. Milton Stewart School of Industrial \& Systems Engineering, Georgia Institute of Technology, Atlanta, GA 30332.
		(email: {\tt gkotsalis3@gatech.edu,   \tt george.lan@isye.gatech.edu}, \tt arkadi.nemirovski@isye.gatech.edu).}
}

\begin{document}

\maketitle

\setcounter{equation}{0}

\begin{abstract}
	
	This work addresses the finite-horizon  robust covariance control problem for discrete-time, partially observable, linear system
	affected  by random zero mean noise and deterministic but unknown disturbances restricted to lie in what is called ellitopic uncertainty set (e.g., finite intersection of centered at the origin ellipsoids/elliptic cylinders).
	Performance specifications are imposed on the random state-control trajectory
	via averaged convex quadratic inequalities, linear inequalities on the mean, as well as  pre-specified upper bounds on the covariance matrix.
	For this problem we develop a computationally tractable procedure for designing affine control policies,  in the sense that the parameters of the policy that guarantees the aforementioned performance specifications are obtained  as solutions to an explicit convex program.
	Our theoretical findings are illustrated by a numerical example.

\end{abstract}

 \vspace{.1in}


 \vspace{.07in}


\vspace{0.1cm}

\section{Introduction}

The standard finite-horizon covariance control problem pertains to designing a feedback control policy for a linear dynamical system affected by zero mean Gaussian disturbances  that will steer
an initial Gaussian random vector (close to) to a prescribed terminal one.  Such stochastic constraints aim to reduce
conservativeness by requiring from the feedback law the imposition of a desired probability distribution profile on the state-control trajectory.
The problem has been studied for linear dynamical systems under full-state feedback   in continuous
\cite{chen_1}, \cite{chen_2}, \cite{chen_3}, as well as discrete-time  settings \cite{bakolas}, \cite{okamoto}.

In this work we consider a robust version of the finite horizon discrete-time covariance control problem under partial observation.
We are given a
discrete-time, linear system,  referred to as $(\mathcal{S})$, with state space description
\begin{eqnarray*}
	\label{system_informal}
	\mathbf{x}_0 & = & \mathbf{z} ~  + ~ \mathbf{s}_0 \\
	\mathbf{x}_{t+1} & = &  A_t  ~ \mathbf{x}_t ~ + ~ B_t ~  \mathbf{u}_t ~ + ~ B_t^{(d)}  ~   \mathbf{d}_t ~+~
	B_t^{(s)} ~ \mathbf{e}_t,~~~~~ \\
	\nonumber
	\mathbf{y}_t  & = & C_t ~  \mathbf{x}_t ~ + ~
	D_t^{(d)} ~  \mathbf{d}_t ~  + ~ D_t^{(s)}   ~  \mathbf{e}_t
	, ~~~t =0,1,2, \hdots,N,
\end{eqnarray*}
where at time $t$, $ x_t \in \bbR^{n_x} $ is the state, $ \mathbf{u}_t \in \bbR^{n_u}$  the control input,
$ \mathbf{y}_t \in \bbR^{n_y}$  the observable output, $ \mathbf{d}_t \in \bbR^{n_d}$ the bounded but unknown exogenous disturbance, while $ \mathbf{e}_t \in \bbR^{n_e}$ is the stochastic exogenous disturbance.
The matrices appearing in the state space description of $(\mathcal{S})$ are known and have compatible dimensions.
In the sequel we treat the initial state $\mathbf{x}_0$ as an exogenous factor with deterministic component $\mathbf{z}$  and a stochastic one $ \mathbf{s}_0 $.
The stochastic exogenous factors form a {\sl random disturbance}  vector
$$
\bfeps=[\mathbf{s}_0;\mathbf{e}_0;...;\mathbf{e}_{N-1}]\qquad \footnotemark \footnotetext{Here and in what follows we use ``MATLAB notation:'' expressions like $z=[u_1;u_2;...;u_k]$ ($z=[u_1,u_2,...,u_k]$) mean that vector/matrix $z$ is obtained from blocks $u_i$ by writing $u_2$ beneath $u_1$, $u_3$ beneath $u_2$, etc. (resp. writing $u_2$ to the right of $u_1$, $u_3$ to the right of $u_2$, etc.)}
$$
and is assumed to have zero mean and known covariance matrix $\Pi$. The deterministic exogenous factors form
the {\sl deterministic disturbance} vector
$$
{\zeta} =[\mathbf{z};\mathbf{d}_0;...;\mathbf{d}_{N-1}]
$$
which is ``unknown but bounded'' -- known to run through a given uncertainty set  $\mathscr{Z}$.

Our goal is to achieve a desired system behavior on the finite horizon $ t = 0, 1, \hdots, N$,  by designing an appropriate non-anticipative affine control law of the form
$$
\mathbf{u}_t = g_{t} + \sum_{i=0}^t
G_i^t~ \mathbf{y}_i,~~
~t =0,1, \hdots, N-1.
$$
where $g_t$ and $ G_{i}^t$ are  vectors and matrices  of appropriate sizes parameterizing the control policy. Performance requirements on $(\mathcal{S})$ are expressed in terms of the resulting random state-control trajectory
$$
\mathbf{w} = [  \mathbf{x}_1;...;\mathbf{x}_{N} ;\mathbf{u}_0;...; \mathbf{u}_{N-1}].
$$
Specifically the desired behavior is modeled by a system of convex constraints on the covariance matrix
$ \Sigma_{\mathbf{w}}$
of the trajectory
of the form
$$
\forall   {\zeta} \in \mathscr{Z}, ~~~~  \mathcal{Q}_i  ~ \Sigma_\mathbf{w}  ~  \mathcal{Q}_i^\T  \preceq \Lambda_i, ~~~~ i = 1 , \hdots, I_c,
$$
where the parameters $ \{ \Lambda_i, \mathcal{Q}_i \}_{i \in \{1, \hdots, I_c\}}$ of compatible dimensions are given.
We will also consider averaged convex quadratic inequality constraints on the state-control trajectory
$$
\forall   {\zeta} \in \mathscr{Z}, ~~~~
\mathbb{E}[\langle \mathcal{A}_i ~ ( \mathbf{w} - \beta_i),    ( \mathbf{w} - \beta_i) \rangle  ~ \leq~ \gamma_i, ~~~~ i = 1 , \hdots, I_q,
$$
where again the parameters $ \{ \mathcal{A}_i, \beta_i, \gamma_i \}_{i \in \{1, \hdots, I_q\} }$ defining the inequalities are assumed to be part of the problem specification.

These two systems of inequalities allow the designer to express among other requirements  also
that the mean $\mu_\mathbf{w} $  of the state-control trajectory is  confined in a pre-specified set, while the covariance $\Sigma_\mathbf{w}$ does  not to exceed a given upper bound, for all realization of the deterministic disturbances $ {\zeta} \in \mathscr{Z}$.

In this paper we develop under
some assumptions on $\mathscr{Z}$, that are for instance satisfied when $\mathscr{Z}$ is the intersection of $K$
concentric ellipsoids/elliptic cylinders,
a computationally tractable procedure for designing
affine control policies, in the sense that the parameters of the policy that guarantees the aforementioned performance
specifications at a certain level of uncertainty are obtained as solutions to an explicit convex program.

Our work expands on
the previous investigations on finite horizon covariance control \cite{chen_1}-\cite{okamoto}
by addressing the robustness issue and the possibility that the full
state may not be available, therefore enabling the steering of the
state-control trajectory density in the presence of disturbances
under partial observation.

Our contribution to the robust control literature lies in the fact that we provide a
tractable computational framework to address the control design problem for a wide class of uncertainty sets
$\mathscr{Z}$ called ellitopes, a concept that was introduced in \cite{jud_17}.
Many interesting uncertainty models
fall in the class of ellitopes, including intersection of ellipsoids and elliptic cylinders, box constraints, as well as the case when the  origin of the exogenous disturbances are discretized continuous-time signals that satisfy Sobolev-type smoothness constraints.
In terms of a tractable formulation we expand considerably on the
 the standard finite-horizon robust control problem utilizing
the induced $l_2$ gain as a sensitivity measure, where the uncertainty set consists of a single ellipsoid (see, e.g., \cite{zhou}, \cite{basar}, \cite{savkin}  and references therein).

From a historical perspective the covariance control problem  in the infinite-horizon case goes back to the original papers on the steady state  assignment problem \cite{hotz1}, \cite{hotz2}, see also \cite{skelton} and the references therein.
The investigation of the finite horizon case was initiated in  \cite{chen_1},  where the optimal covariance steering problem was considered under the assumption that
both input and  noise channels are identical; this requirement was relaxed in \cite{chen_2}.
The work in \cite{bakolas} addressed  the optimal covariance control problem with constraints on input, while \cite{okamoto} included chance constraints in the formulation as well.

This paper considers the finite horizon covariance control problem under partial observation while focusing on the issue of robustness.
Essential to our developments is a key result derived in \cite{jud_17},
specifically, the quantification of the approximation ratio of the semidefinite
relaxation bound on the maximum of a quadratic form over an ellitope.

Another key ingredient to our work is the consideration
of affine policies that are functions of
the so-called ``purified outputs''. This re-parametrization of the set of policies induces a bi-affine structure in the state and control variables that can further be  exploited via robust optimization techniques.
While the specifics are different, in terms of its effect, this
step of passing to purified outputs is akin to the $Q$-parametrization \cite{youla1},\cite{youla2}  in infinite horizon control problems.

Control design in terms of purified outputs was introduced in  \cite{bental2}, where the authors investigated
the finite-horizon robust  control problem for  deterministic linear systems affected by set bounded disturbances
and linear constraints imposed on the state-control trajectory.
In this work we replace linear constraints with convex quadratic ones, while still allowing for a rich class of uncertainty sets, namely ellitopes.
The price one pays for this generalization is that the resulting convex formulation for the affine policy design amounts to a semidefinite program rather than linear programming linear/convex quadratic  programming programs yielded by \cite{bental2}.

The appealing feature of tractability when restricting ourselves to affine policies comes at the cost of potential conservatism.
In the scalar case, under the assumption of deterministic dynamics, the optimality of affine policies for a class of robust multistage
optimization problems subject to constraints on the control variables was proven in
\cite{bertsimas2}.  The degree of suboptimality of affine policies was investigated computationally in \cite{bertsimas3} and \cite{kuhn}. The former reference employs
the sum of squares relaxation, while the latter efficiently computes the suboptimality gap by invoking linear decision rules in the dual formulation as well.
Further applications of affine policies can be found in stochastic programming settings \cite{sun}, \cite{nemirovski}
and in the area of robust model predictive and receding horizon control \cite{loefberg}, \cite{goulart1}, \cite{goulart2},
\cite{skaf}.

The paper is organized as follows.  In Section 2 we introduce the formal problem statement including the notion of ellitopic uncertainty sets, which is new to the control literature. The theoretical developments that lead to a computationally tractable formulation are described
in Section 3, while Section 4 contains an extension to the even wider family of uncertainty sets namely spectratopes. The paper concludes with a numerical example in Section 5. Some technical proofs are relegated to the Appendix.

\subsection{Notation}
The set of real numbers is denoted by $ \bbR$. All vectors are column vectors.
The transpose of the column vector $ x \in \bbR^n$ is written as $ x^\T $ and $ [x]_i$ refers to the $i$-th entry of $ x$, while
$[A]_{i,j}$ refers to the $(i,j)$-th entry of the matrix $ A \in \bbR^{n \times m}$.

The space of symmetric $n\times n$ matrices  is denoted $\mathbb{S}^n$, the {cone of positive semidefinite matrices from $\mathbb{S}^n$  is denoted} $ \mathbb{S}^n_+$, {its interior} $\mathbb{S}_{++}^n$, and the {relations
$A \succeq B$, $A \succ B$ stand for $ A - B \in   \mathbb{S}^n_+$, resp., $A-B\in \mathbb{S}^n_{++}$}.
For $ x, y  \in \bbR^n$, $ \langle x, y \rangle$ denotes the inner product of the two vectors.
A sum (product) of $n  \times n$ matrices over an empty set of indices is the
$n \times n$ zero (identity) matrix.

\section{Problem Statement}

In this section we elaborate on the problem statement described in the introduction, by expanding on the concept of affine policies in purified outputs, the ellitopic uncertainty sets, as well as give some further justification on the use of averaged quadratic inequalities on random state-control trajectory as part of the performance specifications.

\subsection{System dynamics and affine control laws}
We consider a
discrete-time, linear system,  referred to as $(\mathcal{S})$,  with state space description
\begin{eqnarray}
\label{system}
\nonumber
\mathbf{x}_0 & = & \mathbf{z} ~ + ~ \mathbf{s}_0, \\
\mathbf{x}_{t+1} & = &  A_t  ~ \mathbf{x}_t ~ + ~ B_t ~  \mathbf{u}_t ~ + ~ B_t^{(d)}  ~   \mathbf{d}_t ~+~
B_t^{(s)} ~ \mathbf{e}_t \\
\nonumber
\mathbf{y}_t  & = & C_t ~  \mathbf{x}_t ~ + ~
D_t^{(d)} ~  \mathbf{d}_t ~  + ~ D_t^{(s)}   ~  \mathbf{e}_t, ~~~0\leq t \leq N-1.
\end{eqnarray}
where  $x_t\in\bbR^{n_x}$ are states,  $u_t\in\bbR^{n_u}$ are controls,  $y_t\in\bbR^{n_y}$ are observable outputs. The system is affected by two kinds of exogenous factors:\\
-- ``uncertain but bounded'' {\sl deterministic factors} $\mathbf{z}\in\bbR^{n_x}$, $\mathbf{d}_t\in\bbR^{n_d}$, $0\leq t\leq N-1$ which we assemble into a single {\sl deterministic disturbance} vector
$$
\bfzeta=[\mathbf{z};\mathbf{d}_0;...;\mathbf{d}_{N-1}]\in\bbR^{n_\bfzeta},\,n_\bfzeta=n_x+n_dN;
$$
-- {\sl stochastic factors}  $\mathbf{s}_0\in\bbR^{n_x}$, $\mathbf{e}_t\in\bbR^{n_e}$, $0\leq t\leq N-1$, which we assemble into a single {\sl stochastic disturbance} vector
$$
\bfeps=[\mathbf{s}_0;\mathbf{e}_0;...;\mathbf{e}_{N-1}]\in\bbR^{n_\bfeps},\,\,n_\bfeps=n_x+n_eN.
$$
We assume that
\begin{itemize}
\item the matrices $A_t,B_t,...,D_t^{(s)}$ are known for all $t$, $0\leq t\leq N-1$,
\item  deterministic disturbance $\bfzeta$ takes its values in a given {\sl uncertainty set} $\mathscr{Z}\subset\bbR^{n_x+n_dN}$ (of structure to be specified later),
\item  the stochastic disturbance $\bfeps$
is random vector with zero mean and known in advance covariance matrix $\Pi$.
\end{itemize}

Performance specifications will be expressed via constraints on the random state-control trajectory $\mathbf{w}$ and its first two moments. The trajectory is obtained by stacking the
state and the control vectors of various time instances:
\begin{eqnarray}
\label{state-control-trajectory}
\mathbf{x} = [\mathbf{x}_1;...;\mathbf{x}_{N}]\in\bbR^{n_xN},~
\mathbf{u} =  [\mathbf{u}_0;...;\mathbf{u}_{N-1}]\in\bbR^{n_uN},~
\mathbf{w} = [\mathbf{x};\mathbf{u}]\in\bbR^{n_w},n_w=(n_x+n_u)N\
\end{eqnarray}

We will restrict ourselves to casual (non-anticipative)  affine control laws.
The classical output-based (OB) affine policy is
\begin{equation}
\label{control_out}
\mathbf{u}_t = g_{t} + \sum_{i=0}^t
G_i^t~ \mathbf{y}_i,~~
~0\leq t\leq N-1;
\end{equation}
such a policy is specified by the collection of parameters
\begin{equation}
\label{parameters_out_control}
\theta = \{ g_{t}, G_j^t \}_{j \leq t, t \leq N-1},~  g_{t} \in \bbR^{n_u},~ G_j^t \in \bbR^{n_u \times n_x}.
\end{equation}
When employing an (OB) affine policy  with parameters $\theta $
the state-control trajectory is affine
in the initial state and disturbances while highly nonlinear in the parameters $\theta$ of the policy \cite{bental2}.
An essential step towards developing a computationally tractable control design formulation is passing to
deterministic control policies that are affine in the purified outputs of the system, i.e.  purified-output-based (POB) control laws.
To introduce  this concept
 suppose first that the system $(\mathcal{S})$ is controlled by some causal  policy
$$
\mathbf{u}_t = \phi_t( \mathbf{y}_0, \hdots, \mathbf{y}_t), ~~~0\leq t\leq N-1.
$$
We consider now the evolution of a disturbance free system starting from rest, with
state space description
\begin{eqnarray}
\label{model}
\nonumber
\hat{\mathbf{x}}_0 & = & 0, \\
\nonumber
\hat{\mathbf{x}}_{t+1} & = &  A_t ~ \hat{\mathbf{x}}_t ~ + ~ B_t ~ \mathbf{u}_t, \\
\hat{\mathbf{y}}_t  & = & C_t ~ \hat{\mathbf{x}}_t,  ~ 0\leq  t \leq N-1.
\end{eqnarray}
The dynamics in \eqref{model} can be run in ``on-line'' fashion given that
the decision maker
 knows the already computed control values.
The purified outputs $\mathbf{v}_t$ are the differences of the outputs of the actual system \eqref{system} and the outputs of the auxiliary
system \eqref{model}:
\begin{equation}
\label{pure_out}
\mathbf{v}_t = \mathbf{y}_t - \hat{\mathbf{y}}_t , ~  0\leq  t \leq N-1.
\end{equation}
The term ``purified'' stems from the fact that the outputs \eqref{pure_out}
do not depend on the particular choice of control policy.  This circumstance is verified by introducing the signal
$$
{\delta}_t = \mathbf{\mathbf{x}}_t - \hat{\mathbf{x}}_t,~~~0\leq  t \leq N-1,
$$
and observing that the purified outputs can also be obtained
by combining \eqref{system}, \eqref{model}, \eqref{pure_out} and the definition of ${\delta}_t$ as
the outputs of the system with state space description
\begin{eqnarray}
\label{model_pure}
\nonumber
{\delta}_0 & = & \mathbf{z} +\mathbf{s}_0,\\
\nonumber
{\delta}_{t+1} & = & A_t ~ {\delta}_t ~ + ~ B_t^{(d)} ~ \mathbf{d}_t ~ + ~ B_t^{(s)} ~ \mathbf{e}_t,   \\
\mathbf{v}_t  & = & C_t ~ {\delta}_t  ~ + ~ D_t^{(d)} ~ \mathbf{d}_t ~ + ~ D_t^{(s)} ~ \mathbf{e}_t, ~~~0\leq  t \leq N-1.
\end{eqnarray}
Thus, purified outputs are known in advance affine functions of disturbances $\mathbf{\zeta}$, $\mathbf{\epsilon}$.  Furthermore
at time instant $t$,  when the decision on the control value $u_t$ is to be made, we already know the purified outputs $v_0, v_1, \hdots, v_t$.
An affine POB policy is, by definition,  of the form
\begin{equation}
\label{control_pure_out}
\mathbf{u}_t = h_{t} + \sum_{i=0}^t
H_i^t~ \mathbf{v}_i,~~
~0\leq  t \leq N-1;
\end{equation}
such a policy is specified by collection of control parameters
\begin{equation}
\label{parameters_control}
\chi = \{ h_{t}, H_j^t \}_{0\leq j \leq t\leq N-1},~  h_{t} \in \bbR^{n_u},~ H_j^t \in \bbR^{n_u \times n_x}.
\end{equation}
It is proved in \cite{bental2} that {\sl as far as the ``behaviour'' of the controlled system is concerned, the OB and POB control laws are equivalent}. Specifically, equipping  the open loop system $(\mathcal{S})$ with an affine output-based control policy specified by the parameters $\theta$, the state-control trajectory becomes an affine function of the disturbances $[\mathbf{\zeta};\mathbf{\epsilon}]$:
$$
\mathbf{w}=\mathbf{X}_{\theta}\left(\mathbf{\zeta},\mathbf{\epsilon}\right)
$$
for some vector-valued  affine in its argument $[\mathbf{\zeta};\mathbf{\epsilon}]$ function $\mathbf{X}_\theta(\cdot).$ Equipping $(\mathcal{S})$ with an affine POB control policy specified by the parameters $\chi$, the state-control trajectory becomes an affine function of the disturbances:
$$
\mathbf{w}=\mathbf{W}_{\chi}\left(\mathbf{\zeta},\mathbf{\epsilon}\right)
$$
for some vector-valued affine in its argument
$[\mathbf{\zeta};\mathbf{\epsilon}]$ function  $\mathbf{W}_{\chi}(\cdot)$. It is shown in \cite{bental2} that {\sl for every $\theta$ there exists $\chi$ such that $\mathbf{X}_{\theta}(\cdot)\equiv\mathbf{W}_{\chi}(\cdot)$, and vice versa.} It follows that as far as general type affine output based/purified output based control policies are concerned, and whenever the design specifications are expressed in terms of the state-control trajectory (e.g., by imposing on the expectation of the trajectory w.r.t. $\mathbf{\epsilon}$ inequalities which should be satisfied for all realizations $\mathbf{\zeta}\in\mathscr{Z}$), we lose nothing by restricting ourselves to affine POB controls policies. This is what we intend to do, the reason being that, as is easily seen (and is shown in \cite{bental2}), {\sl the mapping  $(\chi,[\mathbf{\zeta};\mathbf{\epsilon}])\to \mathbf{W}_{\chi}\left([\mathbf{\zeta};\mathbf{\epsilon}]\right)$ is bi-affine: affine in $\chi$ when $[\mathbf{\zeta};\mathbf{\epsilon}]$ is fixed, and affine in $[\mathbf{\zeta};\mathbf{\epsilon}]$ when $\chi$ is fixed}. In other words,
\begin{equation}\label{biaff}
\begin{array}{c}
\bfw=\bfW_\chi(\bfzeta,\bfeps)\equiv \bfZ[\chi][\bfzeta;1]+\bfE[\chi]\bfeps\\
\\
\left[\hbox{$\bfZ[\chi]$ and $\bfE[\chi]$: affine in $\chi$}\right]
\\
\end{array}
\end{equation}
As we shall see, this bi-affinity is crucial in the process of reducing the control design  to a computationally tractable problem.
Note that the affine matrix-valued functions $\bfZ[\chi]$, $\bfE[\chi]$ are readily given by the matrices $A_t,...,D_t^{(s)}$, $0\leq  t\leq N-1$, participating in  (\ref{system}).

\subsection{Ellitopic uncertainty sets}

The issue of computational tractability in the ensuing formulation hinges upon the postulated structure of the uncertainty set. Specifically, we will consider a
family of uncertainty sets $\mathscr{D}[\rho]  $,   a single-parameter family of similar to each other ellitopes\footnotemark\footnotetext{an {\sl ellitope} (notion introduced in the context of linear estimation in \cite{jud_17}) is a set which can be represented according to (\ref{ellitop}) as $\mathscr{D}[1]$.} parameterized by the uncertainty level $\rho \geq 0$:
\begin{eqnarray}
\label{ellitop}
\mathscr{D}[\rho] & = &  \{ x  \in   \bbR^{n_{\zeta}} ~|~
\exists (z \in \bbR^{\bar{n}}, t \in \mathcal{T} ): x = P z , ~ \langle S_k z, z \rangle \leq [t]_k \rho, ~ k =1, \hdots, K \},
\end{eqnarray}
where
\begin{itemize}
	\item $\rho > 0 $,  $ P $ is an $ n_{\zeta} \times \bar{n} $  matrix, and  $ S_k $, $1\leq k\leq K$, are $ \bar{n} \times \bar{n}$ matrices with
$S_k\succeq 0$,
	$\sum_{k} S_k \succ 0 $,
	\item $ \mathcal{T}$ is a nonempty computationally tractable convex compact subset of $ \bbR_+^{K}$
	intersecting
	the interior of $ \bbR_+^{K}$
	and such that $ \mathcal{T} $  is monotone, meaning that if $0\leq t\leq t'  $ and $ t' \in \mathcal{T}$ then  $t \in\mathcal{T}$.
\end{itemize}
Note that $\mZ=\sqrt{\rho}\mathscr{D}[1]$.\par
The particular choice of uncertainty description encompasses many interesting cases even from a purely robust control perspective in the absence of stochastic noises in a unifying framework.
\begin{itemize}
	\item When $K = 1 $, $ \mathcal{T} = [0,1]$ and $ Q_1 \succ 0 $, $ \mathscr{D}[\rho] $ reduces to a family of proportional to each other ellipsoids. This class of uncertainty sets is employed in the standard finite-horizon robust control problem utilizing
	the induced $l_2$ gain as a sensitivity measure.
	\item When $ K \geq 1 $, $  \mathcal{T} =  [0,1]^{K}  $, $ \mathscr{D}[\rho] $ is the linear image, under linear mapping $z\to Pz$, of the intersection
$$
\bigcap_{k\leq K} \{ z:  \langle S_k z, z \rangle \leq  \rho \}
$$
of ellipsoids and elliptic cylinders centered at the origin.
\begin{itemize}
	\item When $U \in \bbR^{k \times n}$, $\rank[U] = n$ with rows $u_k^{\T}$, $k \in [K]$ and
	$S_k = u_k u_k^{\T}$, $   \mathscr{D}[\rho] $ is the symmetric with respect to the origin polytope, $ \{ x=Pz, \| U z \|_{\infty} \leq 1\}$.	
\end{itemize}
When for $ p \geq 2$, $ \mathcal{T} = \{ t \in \bbR^{K}_+ : \sum_k [t]_k^{\frac{p}{2}} \leq 1  \} $ and as in the last example $U \in \bbR^{k \times n}$, $\rank[U] = n$ with rows $u_k^{\T}$, $k \in [K]$ and
$S_k = u_k u_k^{\T}$, $   \mathscr{D}[\rho] $ is the set $ \{ x=Pz, \| U z \|_{p} \leq 1\}$.
	\item As another example of an interesting ellitope, we consider the situation where our signals $x$ are discretizations of functions
	of continuous argument running through a compact  domain, and
	the functions we are interested in are those satisfying a Sobolev-type smoothness
	constraints.
\end{itemize}
Note that the family of ellitopes admits ``calculus:'' nearly all operations preserving convexity and symmetry w.r.t. the origin, like taking finite intersections, direct products, linear images, and inverse images under linear embeddings, as applied to ellitopes, result in ellitopes (for details, see \cite[Section 4.6]{A_A_20}).

\subsection{Performance specifications}
\label{PerfSpec}
\an{}{}{
The linear system $(\mathcal{S})$ is subject to both stochastic noise $\bfeps$ and  set-bounded disturbance $\bfzeta$; assuming the control law to be affine POB (this is the default assumption in all our subsequent considerations),  and its state-control trajectory is given by (\ref{biaff}),  $\chi$ being the collection of parameters of the affine POB control law we use.\par
In the sequel, allow for performance specifications expressed by a finite system of
\begin{enumerate}
\item {\sl Averaged over $\bfeps$ convex quadratic constraints on the trajectory} which should be satisfied robustly w.r.t. $\bfzeta\in\mZ$ -- the constraints
\begin{equation}
\label{specification_averaged_quadratic}
\begin{array}{c}
\forall   {\zeta} \in \mathscr{D}[\rho], ~~~
\bbE_{\bfeps}\left[\langle \mathcal{A}_i ~ ( \bfW_\chi(\bfzeta,\bfeps) - \beta_i),    (\bfW_\chi(\bfzeta,\bfeps)- \beta_i) \rangle\right]  ~ \leq~ \gamma_i, ~~~i \in {\cal I}_1\\
\end{array}
\end{equation} and the parameters
\begin{equation}
\label{parameters_averaged_quadratic}
\mathcal{A}_i \in \bbS_+^{n_w},~ \beta_i \in \bbR^{n_w }, ~\gamma_i \in \bbR,\,\,\,i\in{\cal I}_1
\end{equation}
are part of problem's data.  Here and in what follows $\rho>0$ is a given uncertainty level, and ${\cal I}_1,{\cal I}_2,...$ are non-overlapping finite index sets.\\
Note that we do not require for the quadratic forms $\langle\mathcal{A}_i(\bfw-\beta_i),\bfw-\beta_i\rangle$
to exhibit time separability, i.e. stage additivity; thus, our methodology
 is applicable to situations that may not be amenable to dynamic programming.\par
\item {\sl ``Steering of density'' constraints}, specifically,
\begin{enumerate}
\item linear constraints  on the expectation, w.r.t. $\bfeps$, of the state-control trajectory which should be satisfied robustly w.r.t. $\bfzeta\in\mZ$ -- the constraints
\begin{equation}\label{linear}
\begin{array}{c}
\forall {\zeta}\in  \mathscr{D}[\rho]: \langle a_i,\mu_{\bfw}\rangle \leq \gamma_i, \,\,\,i\in{\cal I}_2,\\
\\
\left[\mu_{\bfw}=\mu_{\bfw}[\chi,\bfzeta]:=\bbE_{\bfeps}\left[\bfW_\chi(\bfzeta,\bfeps)\right]\right]\\
\end{array}
\end{equation}
where vectors $a_i$ and scalars $\gamma_i$ are part of the data.
\item convex quadratic constraints on the above expectation $\mu_{\bfw}=\mu_{\bfw}[\chi,\bfzeta]$ which should be satisfied robustly w.r.t. $\bfzeta\in\mZ$ -- the constraints
\begin{equation}
\label{specification_mean_quadratic}
\begin{array}{c}
\forall   {\zeta} \in \mathscr{D}[\rho], ~ \langle \mathcal{A}_i ~  (\mu_\bfw-\beta_i ),   \mu_\bfw-\beta_i  \rangle  ~ \leq~ \gamma_i, ~~~i \in {\cal I}_3\\
\end{array}
\end{equation}
where the parameters
\begin{equation}
\mathcal{A}_i \in \bbS_+^{n_w},~ \beta_i \in \bbR^{n_w }, ~\gamma_i \in \bbR,\,i\in {\cal I}_3,
\end{equation}
are part of the data;
\item $\succeq$-upper bounds on the covariance, w.r.t. $\bfeps$, matrix of the trajectory.  By (\ref{biaff}), {\sl the covariance, w.r.t. $\bfeps$, matrix
$$\Sigma_{\bfw}:=
\mathbb{E}_{\bfeps}\left[(\bfW_\chi(\bfzeta,\bfeps)-\mu_{\bfw}[\chi,\bfzeta])(\bfW_\chi(\bfzeta,\bfeps)-\mu_{\bfw}[\chi,\bfzeta])^\T\right]
$$
of $\bfw=\bfW_\chi(\bfzeta,\bfeps)$ is independent of
$\bfzeta$ and is
\begin{equation}\label{indep}
\Sigma_\bfw=\Sigma_\bfw[\chi]=\bfE[\chi]\Pi\bfE^\T [\chi].
\end{equation}}
$\succeq$-upper bounds on the covariance allowed by our approach are the constraints on $\chi$ of the form
\begin{equation}
\label{specification_covariance}
\begin{array}{c}
\mathcal{Q}_i  ~ \Sigma_\bfw[\chi]  ~  \mathcal{Q}_i^\T  \preceq \Sigma_i,\,\,i\in {\cal I}_4,\\
\end{array}
\end{equation}
where the parameters
\begin{equation}
\Sigma_i \in \bbS^{n_k}_{++},~~
\mathcal{Q}_i \in \bbR^{n_k \times n_w }\,\,i\in{\cal I}_4,
\end{equation}
are part of problem's data.
\end{enumerate}
\end{enumerate}
A particular scenario where constraints (\ref{specification_averaged_quadratic}) are relevant is when the decision maker has some desired targets for the states and controls along with allowed deviation levels
to be satisfied robustly for all $ {\zeta} \in \mathscr{D}[\rho]$; to handle these specifications, it suffices to include in the list of quadratic functions of $\bfw$ participating in (\ref{specification_averaged_quadratic}) the functions
$$
	\mathbb{E}_{\bfeps}[\| \mathbf{x}_t - \hat{x}_t \|_2^2]  \leq \hat{\gamma}_t,\,\,
\mathbb{E}_{\bfeps}[\| \mathbf{u}_t - \tilde{u}_t \|_2^2]  \leq \tilde{\gamma}_t
$$
with relevant $t$'s. Similarly, constraints (\ref{specification_mean_quadratic}) allow to handle robust, w.r.t. $\bfzeta\in\mZ$, versions of the constraints
$$
	\| \mu_{\mathbf{x},t} - \hat{\mu}_t \|^2   \leq   \hat{\epsilon}_t,  \,
	\| \mu_{\mathbf{u},t} - \tilde{\mu}_t \|^2  \leq    \tilde{\epsilon}_t,
$$
where $\mu_{\mathbf{x},t}=\mu_{\mathbf{x},t}[\chi,\bfzeta]$ is the expectation, w.r.t. $\bfeps$, of the $t$-th state $x_t$, and similarly for $\mu_{\mathbf{u},t}$. Finally,
constraints \eqref{specification_covariance} with appropriately chosen weighting matrices $\{ \mathcal{Q}_i \}_{i\in{\cal I}}$ allow to impose $\succeq$-upper bounds
on the covariance, w.r.t. $\bfeps$, matrices of states and/or controls at prescribed time instants.
\begin{eqnarray*}
	\Sigma_{\mathbf{x},t}  \preceq \hat{\Lambda}_t, ~ t \in T_x,~
	\Sigma_{\mathbf{u},t} \preceq \tilde{\Lambda}_t, ~ \forall t \in T_u,
\end{eqnarray*}
where $ \{ \hat{\Lambda}_t \}_{t \in T_x}$,  $ \{ \tilde{\Lambda}_t \}_{t \in T_u}$ are prespecified positive definite matrices at  particular instants indexed by $T_x$ and $T_u$ respectively.
\paragraph{The problem} we are interested in is
\begin{equation}\label{probofint}
\hbox{\sl Find  affine  POB $\chi$ satisfying a given set of constraints {\rm  (\ref{specification_averaged_quadratic}),   (\ref{linear}), (\ref{specification_mean_quadratic}), (\ref{specification_covariance})}}.
\end{equation}
As stated, our problem of interest is a feasibility problem; the approach we are about to develop can be straightforwardly applied to the optimization version of our problem, where the right hand sides in the constraints just specified, instead of being part of the data, are treated as additional variables, and our goal is to minimize a convex function of $\chi$ and these additional variables under the constraints stemming from
 (\ref{specification_averaged_quadratic}),  (\ref{linear}),  (\ref{specification_mean_quadratic}),  (\ref{specification_covariance}) (and, perhaps, additional convex constraints linking $\chi$ and the added variables).

\section{Convex design}

\subsection{Processing performance specifications}\label{processingspec}
We are about to  process one by one the performance specifications indicated in Section \ref{PerfSpec}.
\begin{enumerate}
\item {[constraints  (\ref{specification_averaged_quadratic})]} Given matrix $\mathcal{A}\succeq0$ and a vector $\beta$ and taking into account(\ref{biaff}) and the fact that $\bfeps$  is zero mean with covariance matrix $\Pi$, we have
\begin{equation}\label{have1}
\begin{array}{l}
\mathbb{E}_{\bfeps}\left[\langle \mathcal{A}(\bfW_{\chi}(\bfzeta,\bfeps)-\beta),(\bfW_{\chi}(\bfzeta,\bfeps)-\beta)\rangle\right]\\
\equiv
[\bfzeta;1]^\T  \bfZ^\T [\chi]\mathcal{A}\bfZ[\chi][\bfzeta;1]-2\beta^\T\cA\bfZ[\chi][\bfzeta;1]+\Tr(\Pi^{1/2}\bfE^\T [\chi]\mathcal{A}\bfE[\chi]\Pi^{1/2})+r\\
\end{array}
\end{equation}
with $r$ readily given by $\mathcal{A}$ and $\beta$.
As a result (for verification, see Section \ref{verifyA})
\begin{quote}
{\bf A.} {\sl Given matrix $\mathcal{A}\succeq0$, vector $\beta$, and scalar $\gamma$, the semi-infinite constraint
\begin{equation}\label{constraintA}
\forall \bfzeta\in\mZ: \mathbb{E}_{\bfeps}\left[\langle \mathcal{A}(\bfW_{\chi}(\bfzeta,\bfeps)-\beta),(\bfW_{\chi}(\bfzeta,\bfeps)-\beta )\rangle\right]\leq\gamma
\end{equation}
in variables $\chi$ is equivalently represented by the system of constraints
\begin{equation}\label{systemA}
\begin{array}{cl}
\left[\begin{array}{c|c}X&x\cr\hline x^\T &\alpha\cr\end{array}\right]\succeq \bfZ^\T [\chi]\mathcal{A}\bfZ[\chi]&(a)\\
\\
\alpha+\delta+\|\mathcal{A}^{1/2}\bfE[\chi]R\|_{\Fro}^2+r\leq \gamma&(b) \\
\\
\forall (\bfzeta\in\mZ): \bfzeta^T X\bfzeta+2\bfzeta^\T [x+p[\chi]]\leq\delta-q[\chi]&(c)\\
\end{array}
\end{equation}
in variables $\chi,X,x,\alpha,\delta$; here $\|\cdot\|_{\Fro}$ is the Frobenius norm of a matrix, $R=\Pi^{1/2}$, and
$
p[\chi]$, $q[\chi]$ are affine in $\chi$ vector-and real-valued functions uniquely defined by the identity
\begin{equation}\label{ident1}
\forall z\in\bbR^{n_\zeta}: -2\beta^\T\cA\bfZ[\chi][z;1]=2p^\T[\chi]z+q[\chi].
\end{equation}
Note that
$r$ and the functions $p[\cdot]$, $q[\cdot]$ are readily given by $\cA$, $\beta$ and problem's data.}
\end{quote}

Here and in what follows, claim of the form  ``system $(P)$ of constraints on variables $x$ is equivalently represented by system $(R)$ of constraints on variables $x,y$'' means that  a candidate solution $x$   is feasible for $(P)$ {\sl if and only if} $x$ can be extended, by appropriately selected $y$, to a feasible solution of $(R)$. Whenever this is the case, finding a feasible solution to $(P)$ straightforwardly reduces to finding a feasible solution to $(R)$, and similarly for optimizing an objective $f(x)$ over the feasible solutions of $(P)$.

\item {[constraints  (\ref{linear})]} By (\ref{biaff}), we have $\mu_{\bfw}=\mu_{\bfw}[\chi,\bfzeta]=\bfZ[\chi][\bfzeta;1]$. As a result, for a vector $a$ and real $\gamma$ we have
\begin{equation}\label{have2}
\begin{array}{lll}
&\forall (\bfzeta\in\mZ): \langle a,\mu_{\bfw}\rangle \leq\gamma&(a)\\
\\
\Leftrightarrow&\forall (\bfzeta\in\mZ): a^\T\bfZ[\chi][\bfzeta;1] \leq\gamma&(b)\\
\end{array}
\end{equation}
As a result (for verification, see Section \ref{verifyAA})
\begin{quote}
{\bf B.} {\sl Given vector $a$ and scalar $\gamma$, the semi-infinite constraint
\begin{equation}\label{constraintAA}
\forall \bfzeta\in\mZ: \langle a,\mu_{\bfw}[\chi,\bfzeta]\rangle \leq\gamma
\end{equation}
in variables $\chi$ is equivalently represented by the system of convex constraints
\begin{equation}\label{systemAA}
\begin{array}{cl}
\lambda=[\lambda_1;...;\lambda_K]\geq0&(a)\\
\\
\left[\begin{array}{c|c}\alpha&-p^\T[\chi]P\cr\hline -P^\T p[\chi]&\sum_k\lambda_kS_k\cr\end{array}\right]\succeq0&(b)\\
\\
\alpha+\rho\phi_{\cT}(\lambda)+q[\chi]\leq\gamma&(c)\\
\end{array}
\end{equation}
in variables $\chi,\lambda,\alpha$; here $
p[\chi]$, $q[\chi]$ are affine in $\chi$ vector-valued and real-valued functions uniquely defined by the identity
\begin{equation}\label{ident2}
\forall z\in\bbR^{n_\zeta}: a^\T\bfZ[\chi][z;1]=2z^\T p[\chi]+q[\chi],
\end{equation}
$\phi_\cT(\lambda)$ is the support function of $\cT$:
$$
\phi_\cT(\lambda)=\max_{t\in\cT} \lambda^\T t,
$$
and $P$ and $S_k$ come from (\ref{ellitop}).
Note that the functions $p[\cdot]$, $q[\cdot]$ are readily given by $a$ and problem's data.}
\end{quote}

\item {[constraints (\ref{specification_mean_quadratic})]} As it was already mentioned, $\mu_{\bfw}=\bfZ[\chi][\bfzeta;1]$. As a result, for a matrix $\mathcal{A}\succeq0$ and a vector $\beta$ we have
$$
\begin{array}{l}
\langle \mathcal{A}(\mu_{\bfw}[\chi,\bfzeta]-\beta),\mu_{\bfw}[\chi,\bfzeta]-\beta\rangle
\equiv\langle \mathcal{A}(\bfZ[\chi][\bfzeta;1]-\beta), \bfZ[\chi][\bfzeta;1]-\beta)\\
\equiv
[\bfzeta;1]^\T  \bfZ^\T [\chi]\mathcal{A}\bfZ[\chi][\bfzeta;1]-2\beta^\T\mathcal{A}\bfZ[\chi][\bfzeta;1]+r
\end{array}
$$
with scalar $r$ readily given by $\cA$, $\beta$. Applying the reasoning completely similar to the one used to justify claim {\bf A}, see Section \ref{verifyA},
we conclude that
\begin{quote}
{\bf C.} {\sl Given matrix $\mathcal{A}\succeq0$ and  vector $\beta$, the semi-infinite constraint
\begin{equation}\label{constraintB}
\forall \bfzeta\in\mZ: \langle \mathcal{A}(\mu_{\bfw}[\chi,\bfzeta]-\beta),\mu_{\bfw}[\chi,\bfzeta]-\beta\rangle\leq\gamma
\end{equation}
in variables $\chi$ can be equivalently represented by the system of constraints
\begin{equation}\label{systemB}
\begin{array}{cl}
\left[\begin{array}{c|c}X&x\cr\hline x^\T &\alpha\cr\end{array}\right]\succeq \bfZ^\T [\chi]\mathcal{A}\bfZ[\chi]&(a)\\
\\
\alpha+\delta+r\leq\gamma&(b)\\
\\
\forall (\bfzeta\in\mZ): \bfzeta^\T  X\bfzeta+2\bfzeta^\T [x+p[\chi]]\leq \delta-q[\chi]&(c)\\
\end{array}
\end{equation}
in variables $\chi,X,x,\alpha,\delta$, with and affine vector- and real-valued functions $p[\cdot]$ and $q[\cdot]$ given by the identity (\ref{ident1}).}
\end{quote}
Same as in {\bf A} and {\bf B}, scalar $r$ and functions $p[\cdot]$, $q[\cdot]$ are readily given by $\cA$, $\beta$ and problem's data.
\item {[constraints (\ref{specification_covariance})]} Finally, from (\ref{indep}) and (\ref{specification_covariance}) it follows that
\begin{quote}
{\bf D.} {\sl An $\succeq$-upper bound on the state-control trajectory's covariance matrix of the form
$$
\mathcal{Q}\Sigma_{\bfw}[\chi]\mathcal{Q}^\T \preceq\Sigma
$$
or, equivalently,
$$
[\mathcal{Q}\bfE[\chi]R][\mathcal{Q} \bfE[\chi]R]^\T \preceq\Sigma, \eqno{[R=\Pi^{1/2}]}
$$
see (\ref{indep}), is nothing but the Linear Matrix Inequality
\begin{equation}\label{systemC}
\left[\begin{array}{c|c}\Sigma&\mathcal{Q}\bfE[\chi]R\cr\hline
R^\T \bfE^\T [\chi]\mathcal{Q}^\T &I_{n_{\bfeps}}\cr\end{array}\right]\succeq0
\end{equation}
(we have used the Schur Complement Lemma).}
\end{quote}
\end{enumerate}
\paragraph{Intermediate conclusion.} We arrive at the following corollary of {\bf A} -- {\bf D}:
\begin{corollary}\label{intcor} The problem of interest can be equivalently reformulated
as an explicit convex feasibility problem $(\mathcal{F})$ in variables $\chi$ specifying the control law we are looking for and in several additional variables. The constraints of $(\mathcal{F})$ can be split into two groups:
\begin{itemize}
\item ``simple'' --- explicit computationally tractable convex constraints like {\rm (\ref{systemAA}), (\ref{systemC})} or constraints $(a)$, $(b)$ in systems {\rm (\ref{systemA}), (\ref{systemB})};
\item ``complex'' -- semi-infinite constraints $(c)$ coming from  {\rm (\ref{systemA}), (\ref{systemB})}.
\end{itemize}
\end{corollary}

\noindent If there were no constraints we have called ``complex,'' the situation would be extremely simple --- the feasibility problem of interest (and its optimization versions, where the right hand sides of the constraints
(\ref{specification_averaged_quadratic}), (\ref{linear}), (\ref{specification_mean_quadratic}), (\ref{specification_covariance})) are, same as $\chi$, treated as design variables and we want to minimize convex function of
all these variables under the constraints stemming from (\ref{specification_averaged_quadratic}), (\ref{linear}), (\ref{specification_mean_quadratic}), (\ref{specification_covariance}), and, perhaps, additional explicit convex constraints) would be just an explicit and therefore efficiently solvable convex optimization problem.
We are about to demonstrate that
\begin{enumerate}
\item In the {\sl simple case}, when $\mZ$ is an ellipsoid (i.e., $K=1$ in (\ref{ellitop})), the complex constraints admit simple equivalent representation;
\item Beyond the simple case, complex constraints can be computationally intractable, but they admit reasonably tight ``safe'' tractable approximations, resulting in computationally efficient ``moderately conservative'' design of affine POB policy.
\end{enumerate}

\subsection{Processing complex constraints}
In this Section, our goal is to process a semi-infinite constraint
$$
\forall (\bfzeta\in\mZ): \bfzeta^\T  X\bfzeta+2x^\T \bfzeta\leq \xi
\eqno{(C[\rho])}
$$
in variables $X,x,\xi$. It is easily seen that all complex constraints in the system $(\cF)$ we are interested in (in our context, these are constraints $(c)$ in (\ref{systemA}), (\ref{systemB})) are of the generic form $(C[\rho])$, with $X,x,\xi$ being affine functions of the variables of $(\cF)$.

\subsubsection{Simple case: $K=1$ in (\ref{ellitop})}
In the simple case $K=1$, we have
$$
\mZ=\{\bfzeta=Pz: z^\T  \bar{Q}z\leq 1\},
$$
where $\bar{Q}\succ0$ (in terms of (\ref{ellitop}), $\bar{S}=Q_1/(\rho\max_{t\in\mathcal{T}}t_1)$).
By Inhomogeneous $\cal{S}$-Lemma (see, e.g., \cite{bental6}), $(C[\rho])$ can be equivalently represented by the Linear Matrix Inequality
\begin{equation}\label{caseK=1}
\left[\begin{array}{c|c}P^\T XP+\lambda\bar{Q}&P^\T x\cr\hline
x^\T P&\xi-\lambda\cr\end{array}\right]\succeq0\ \&\ \lambda\geq0
\end{equation}
in variables $X,x,\xi,\lambda$. In view of Corollary \ref{intcor}, we conclude that the problem of interest can be equivalently reformulated as feasibility problem with efficiently computable convex constraints, implying efficient solvability of the problem.

\subsubsection{General case: $K>1$ in (\ref{ellitop})}
It is well known that in the general case, already checking the validity of $(C[\rho])$ for given $X,x,\xi$ can be NP-hard. However, $(C[\rho])$ admits a ``tight,'' in certain precise sense, tractable approximation
given by the following result
\begin{theorem}\label{thean1}
(i) Let, same as above,
$$
\phi_{\mathcal{T}}(\lambda)=\max_{t\in\mathcal{T}}\lambda^\T t
$$ be the support function of $\cT$.
The system of convex constraints
$$
\begin{array}{cl}
\lambda_k\geq0,k\leq K&(a)\\
\\
\left[\begin{array}{c|c}\sum_k\lambda_kQ_k-X&-x\cr\hline -x^\T  &\tau\cr\end{array}\right]\succeq0&(b)\\
\\
\rho\phi_{\mathcal{T}}(\lambda)+\tau\leq\xi&(c)\\
\end{array}
\eqno{(S[\rho])}
$$
in variables $X,x,\xi,\lambda=\{\lambda_k,k\leq K\},\tau$ is safe tractable approximation of $(C[\rho])$, meaning that whenever  $(X,x,\xi)$ can be augmented by properly selected $\lambda$, $\tau$ to a feasible solution  to $(S[\rho])$,
$X,x,\xi$ is feasible for $(C[\rho])$.
\par(ii) The safe tractable approximation $(S[\rho])$ of $(C[\rho])$ is tight within the factor
\begin{equation}\label{tight}
 \varkappa=3\ln(6K)
\end{equation}
meaning that whenever $(X,x,\xi)$ can\underline{not} be augmented by $\lambda$ and $\tau$ to yield a feasible solution to $(S[\rho])$, $(X,x,\xi)$ is \underline{not} feasible for $(C[\varkappa\rho])$.
\end{theorem}
\par\noindent
For proof, see Appendix \ref{appextSlemmaEll}.
\par
As a consequence of Theorem \ref{thean1} we arrive at the following {\sl conservative} approximation of our problem of interest: we build the system $(\mathcal{S}[\rho])$ of explicit convex constraints
on $\chi$ and additional variables by augmenting the simple, in terms of Corollary \ref{intcor}, constraints \begin{itemize}
\item in simple case $K=1$ --- by equivalent representations (\ref{caseK=1}) of complex, in terms of the same Corollary,  constraints;
\item In the case of $K>1$ --- by safe tractable approximations, as given by Theorem \ref{thean1}, of complex constraints.
\end{itemize}
With this approach, the variables $\chi$ participating in the problem of interest form part of the variables of the efficiently solvable problem $(\mathcal{S}[\rho])$, and
the $\chi$-component of every feasible solution, if any, to the latter problem is feasible for the problem of interest. In the case of $K=1$, we can get in this way {\sl all}
feasible solutions to the problem of interest; in particular, $(\mathcal{S}[\rho])$ is infeasible if and only if the problem of interest is so. In the case of $K>1$, it may happen that
$(\mathcal{S}[\rho])$ is infeasible, while the problem of interest is not so; however, in the latter case the problem of interest becomes infeasible after the original level of uncertainty $\rho$ is increased by ``moderate'' factor, namely, to $\varkappa\rho$ (that is, when the original uncertainty set $\mZ$ is increased by factor $\sqrt{\varkappa}$).
\par
Another, sometimes more instructive, way to express the same fact is as follows. Consider the optimization version ${\cal P}[\rho]$ of the problem of interest, where the right hand sides, let they be called $\gamma$, of the original constraints
become additional to $\chi$ variables, and one seeks to minimize a given efficiently computable convex function of $\gamma$ and $\chi$ under the original constraints
and, perhaps, system $(A)$ of additional efficiently computable convex constraints on $\gamma$ and $\chi$. As is immediately seen, $\gamma$'s are among the right hand sides of the constraints of $(\mathcal{S}[\rho])$, and we can apply similar procedure to the latter system -- treat $\gamma$'s as additional variables  and optimize the same objective over all resulting variables under the constraints stemming from $(\mathcal{S}[\rho])$ and additional constraints $(A)$, if any. With this approach, every feasible solution to the resulting optimization problem, let it be called ${\cal P}^+[\rho]$, induces a feasible solution, with the same value of the objective, to ${\cal P}[\rho]$, and we have
\def\Opt{{\mathrm{Opt}}}
$$
\Opt[{\cal P}[\rho]]\leq\Opt[{\cal P}^+[\rho]]\leq \Opt[{\cal P}[\varkappa\rho]]\,\,\forall \rho>0,
$$ where $\Opt[{\cal P}]\in\bbR\cup\{+\infty\}$ stands for the optimal value of an optimization problem ${\cal P}$.
}

\section{Extension: Spectratopic uncertainty}
The above constructions and results can be straightforwardly extended to the case when the single-parametric family $\mZ$ of uncertainty  sets is comprised of
similar to each other {\sl spectratopes} rather than similar to each other ellitopes, that is,
\begin{equation}\label{spectratop}
\mZ=\{x\in\bbR^{n_\zeta}:\exists (t\in\mathcal{T},z\in\bbR^{\bar{n}}): x=Pz, S_k^2[z]\preceq \rho t_kI_{f_k},k\leq K\},
\end{equation}
where $\mathcal{T}$ is as in (\ref{ellitop}), and $S_k[\cdot]$ are linear functions taking values in the space $\bbS^{f_k}$:
\begin{equation}\label{Ski}
S_k[z]=\sum_{i=1}^{\bar{n}} z_iS^{ki},\,\, S^{ki}\in\bbS^{f_k}.
\end{equation}
Spectratopes -- the sets which can be represented according to (\ref{spectratop}) as $\mathscr{D}[1]$ --- were introduced in \cite{jud_18}. Every ellitope is a spectratope, but not vise versa; an
instructive example of a ``genuine'' spectratope is the matrix box -- the unit ball of the spectral norm $\|\cdot\|_{2,2}$:
$$
\mathcal{B}=\{z\in\bbR^{p\times q}: \|z\|_{2,2}\leq 1\}=\{z\in\bbR^{p\times q}: \exists t_1\in\cT=[0,1]: \underbrace{\left[\begin{array}{c|c}&z\cr\hline z^\T &\cr\end{array}\right]^2}_{S_1^2[z]}\preceq t_1I_{p+q}\}.
$$
Spectratopes allow for the same rich calculus as ellitopes, see \cite[Section 4.6]{A_A_20}.
\par
Theorem \ref{thean1} admits the ``spectratopic extension'' as follows (for proof, see Appendix \ref{appextSlemma}):
\begin{theorem}\label{thean2} Consider the semi-infinite constraint
$(C[\rho])$ with $\mZ$ given by {\rm (\ref{spectratop})} rather than by {\rm (\ref{ellitop})}. Then \par
(i) For $k\leq K$, let the mappings $Z\mapsto{\cal S}_k[Z]:\bbS^{\bar{n}}\to\bbS^{f_k}$ and $D\mapsto {\cal S}_k^*[D]: \bbS^{f_k}\to\bbS^{\bar{n}}$ be given by
$$
{\cal S}_k[Z]=\sum_{i,j}Z_{ij}S^{ki}S^{kj},\,\, [{\cal S}_k^*[D]]_{ij}=\Tr(DS^{ki}S^{kj})
$$
where $S^{ki}$ are given by {\rm (\ref{spectratop}), (\ref{Ski})}. The system of convex constraints in variables $X,x,\xi,\Lambda=\{\Lambda_k,k\leq K\},\tau$
$$
\begin{array}{cl}
\Lambda_k\in\bbS^{f_k}_+,k\leq K&(a)\\
\\
\left[\begin{array}{c|c}\sum_k{\cal S}_k^*[\Lambda_k]-X&-x\cr\hline-x^\T  &\tau\cr\end{array}\right]\succeq0&(b)\\
\\
\rho\phi_{\mathcal{T}}(\lambda[\Lambda])+\tau\leq\xi&(c)\\
\end{array}
\eqno{(\bS[\rho])}
$$
where
$$
\lambda[\Lambda]=[\Tr(\Lambda_1);...;\Tr(\Lambda_{K})]
$$
and $\phi_{\cT}$ is the support function of $\cT$
 is safe tractable approximation of $(C[\rho])$, meaning that whenever  $(X,x,\xi)$ can be augmented by properly selected $\Lambda$, $\tau$ to a feasible solution  to $(S[\rho])$,
$(X,x,\xi)$ is feasible for $(C[\rho])$.
\par(ii) The safe tractable approximation $(\bS[\rho])$ of $(C[\rho])$ is tight within the factor
\begin{equation}\label{tightnew}
 \varkappa=2\ln\left(8{\sum}_{k=1}^{K}f_k\right),
\end{equation}
meaning that whenever $(X,x,\xi)$ can\underline{not} be augmented by $\Lambda$ and $\tau$ to yield a feasible solution to $(\bS[\rho])$, $(X,x,\xi)$ is \underline{not} feasible for $(C[\varkappa\rho])$.
\end{theorem}
\par\noindent
Claim {\bf B} from Section \ref{processingspec} also admits a spectratopic extension, specifically,
\begin{quote}
$\mathbf{B'}$
{\sl Given vector $a$ and scalar $\gamma$, the semi-infinite constraint
\begin{equation}\label{constraintAAA}
\forall \bfzeta\in\mZ: \langle a,\mu_{\bfw}[\chi,\bfzeta]\rangle \leq\gamma
\end{equation}
in variables $\chi$ with $\mZ$ given by (\ref{spectratop}) is equivalently represented by the system of convex constraints
\begin{equation}\label{systemAA_spectr}
\begin{array}{cl}
\Lambda=\{\Lambda_k\in\bbS^{f_k}_+,k\leq K\}&(a)\\
\\
\left[\begin{array}{c|c}\alpha&-p^\T[\chi]P\cr\hline -P^\T p[\chi]&\sum_k\cS^*_k[\Lambda_k]\cr\end{array}\right]\succeq0&(b)\\
\\
\alpha+\rho\phi_{\cT}(\lambda[\Lambda])+q[\chi]\leq\gamma&(c)\\
\end{array}
\end{equation}
in variables $\chi,\lambda,\alpha$; here $\cS_k^*[\cdot]$ and $\lambda[\Lambda]$ are as in Theorem \ref{thean2}, and $
p[\chi]$, $q[\chi]$ are affine in $\chi$ vector-valued and real-valued functions uniquely defined by the identity (\ref{ident2}).
Note that the functions $p[\cdot]$, $q[\cdot]$ are readily given by $a$ and problem's data.}
\end{quote}
(for proof, see Section \ref{verifyBprime}).
\par
The consequences of $\mathbf{B'}$ and Theorem \ref{thean2} for design of affine POB control policies are completely similar to those of Theorem \ref{thean1}, with the scope of these consequences wider than before due
the
significant extension of the family of uncertainty sets we can handle.

\section{Numerical Illustration }

In this section we present a numerical example of the design procedure. We consider the output-feedback control problem
of the longitudinal dynamics of a Boeing 747 aircraft. The continuous-time dynamics are linearized and are valid for modest deviations from the trim conditions of $40000$ ft altitude and $851.3$ ft/sec ground speed, the tailwind  is $73.3$ ft/sec. The aircraft is flying horizontally in the $ XZ$-plane, and its longitudinal axis coincides with the $X$-axis. The data are taken from  \cite{boyd_lecture} and are reported in \cite{bental6} as well.
The state vector  is
$$
x(t)^\T  = \AR{ccccc}{ \Delta h(t) & \Delta u(t) & \Delta v(t) & q (t) &  \Delta \theta(t) },
$$
\begin{itemize}
	\item $\Delta h$ : deviation in altitude [ft], positive is down
	\item $\Delta u$ : deviation in velocity along the aircraft’s axis [ft/sec], positive is forward
	\item $\Delta v$ : deviation in velocity orthogonal to the aircraft's axis [ft/sec], positive is down
		\item $\Delta \theta$  : deviation in angle $\theta$ between the aircraft’s axis and the $X$-axis
	positive is up, [crad]   the unit is 0.01 radian
\item $q$: angular velocity of the aircraft (pitch rate), [crad/sec].
\end{itemize}
The control vector is
$$
u(t)^\T  = \AR{ccccc}{ \delta_e(t) & \delta_{th}(t) },
$$

\begin{itemize}
\item $\delta_e$    : elevator (control surface) angle deviation, positive is down  [crad]
\item $\delta_{th}$ : engine thrust deviation [10000 lbs]
\end{itemize}
The disturbance vector is
$$
d(t)^\T  = \AR{ccccc}{ u_w(t) & v_w(t) },
$$
\begin{itemize}
	\item $u_w$ : wind velocity along the aircraft’s axis, [ft/sec].
	\item $v_w$ : wind velocity orthogonal to the aircraft’s axis, [ft/sec].
\end{itemize}
The output vector is
$$
y(t)^\T  = \AR{ccccc}{ \delta u(t) & \dot{h}(t) },
$$
\begin{itemize}
	\item $\delta u$ : velocity deviation, , [ft/sec]
	\item $\dot{h}$ : climb rate, [ft/sec]
\end{itemize}
The state space description of the continuous time model is
 \begin{eqnarray*}
 	\frac{d}{dt} x(t) & = &
 	A_c ~
 x(t)
 	~ + ~
 	B_c ~
 u(t)~  +  ~ D_c~ d(t),
 	~~~~y(t) =  C_c ~
 	x(t),
 \end{eqnarray*}
with the system matrices
{\small\begin{eqnarray*}
A_c =
\AR{rrrrr}{
0 & 0 & - 1 & 0 & 0 \\
0 & -0.003 & 0.039 & 0 & -0.322 \\
0 & -0.065 & -0.319 & 7.74 & 0 \\
0 & 0.020 & -0.101 & - 0.429 & 0 \\
0 & 0 & 0 & 1 & 0
},~
B_c =  \AR{rr}{
	0 &  0  \\
	0.01 & 1 \\
	-0.18 & -0.04 \\
	-1.16 & 0.598  \\
	0 &  0
},~
D_c =
\AR{rr}{
	0  & - 1  \\
	0.003 & -0.039 \\
	0.065 & 0.319 \\
	-0.020 & 0.101 \\
	0 & 0
}
\end{eqnarray*}
and  observation matrix
$$
C_c = \AR{rrrrr}{
	0 & 1 &  0 & 0 & 0\\
	0 & 0 & -1 & 0 & 7.74
}.
$$
We discretize the continuous-time dynamics at a sampling rate of $T = 10$ sec assuming that all the input signals
are piecewise-constant, i.e.
$$
u(t)  \equiv u_k~~~~~~~~~~\text{for}~~   k T \leq t < (k+1) T.
$$
The discrete-time linear time-invariant state space model is
\begin{eqnarray*}
	\label{Boeing}
	\mathbf{x}_{t+1} & = &  A  ~ \mathbf{x}_t ~ + ~ B ~  \mathbf{u}_t ~ + ~ B^{(d)}  ~   (\mathbf{d}_t  + \mathbf{e}_t ) =  A  ~ \mathbf{x}_t ~ + ~ B ~  \mathbf{u}_t ~ + ~ B^{(d)}  ~   \mathbf{d}_t  ~ + ~ B^{(s)}  \mathbf{e}_t    \\
	\nonumber
	\mathbf{y}_t  & = & C ~  \mathbf{x}_t
	, ~~~t =0,1,2, \hdots,N,
\end{eqnarray*}
with
$$
A = \exp[A_c T],~ B = \int_{0}^\T  \exp[A (T- \sigma)] B_c d \sigma,   ~ B^{(s)} =B^{(d)} = \int_{0}^\T  \exp[A (T- \sigma)] D_c d \sigma,  ~ C = C_c.
$$
The state space matrices compute to be
$$
A =
\AR{rrrrr}{
1  & -1.1267  & -0.6528   & -8.0749 &   1.5890 \\
0  &  0.7741  &  0.3176   & -0.9772 &  -2.9690 \\
0  &  0.1157  &  0.0201   & -0.0005 & -0.3628  \\
0  &  0.0111  &  0.0033   & -0.0349 &  -0.0447 \\
0  &  0.1388  & -0.0862   &  0.2935 &   0.7579 },~
B = \AR{rr}{
89.1973 & -50.1685 \\
5.2231  &   6.3614 \\
-9.4731 &   5.9294 \\
-0.3236 &   0.3178 \\
-4.5318 &   3.2146	},
$$
$$
B^{(d)}  =\AR{rr}{
 1.1267 & -19.3472 \\
 0.2259 &  -0.3176 \\
-0.1157 &   0.9799 \\
-0.0111 &  -0.0033 \\
-0.1388 &   0.0862}.
$$}\noindent
The linearized flight dynamics correspond to the operating conditions of
$ 528$  mph aircraft speed and $ 50 $ mph steady-state tail wind.
The initial state is
$$
\mathbf{x}_0 = [0;0 ;0 ; 0 ;0].
$$
Our  goal is to maintain  the operating conditions at $t = 100$ sec  and $ t = 200 $sec.
We set the horizon at $N=20$, and the target states
$$
\tilde{\mathbf{x}}_{10}   =   \tilde{\mathbf{x}}_{20} = \mathbf{x}_0.
$$
We assume that the stochastic noise is standard Gaussian (zero mean, unit covariance: $\Sigma_{[\mathbf{e}_0; \mathbf{e}_1 ; ... ; \mathbf{e}_{N-1}]} =  I $), while the deterministic disturbances satisfy the constraints
$$
\sum_{t=0}^{N-5}  \|\mathbf{d}_t \|_2^2 \leq 1, ~     \sum_{t=N-4}^{N-1}  \|\mathbf{d}_t \|_2^2   \leq  \frac{1}{100},
$$
reflecting the circumstance that the near horizon exogenous disturbances are more severe.
The performance constraints are
$$
 \bbE[ \| \mathbf{x}_{10} -  \tilde{\mathbf{x}}_{10} \|^2 ] \leq 400, ~  \bbE[ \| \mathbf{x}_{N} -  \tilde{\mathbf{x}}_{N} \|^2 ] \leq 400, ~~ \mathbf{\Sigma}_{\mathbf{x}_N} \preceq 400 I
$$
and are to be met robustly.
We design a POB control law that meets the performance specifications.
We run the computations in a MATLAB environment on
a MacBook Pro with a 2.3 GHz Intel Core i5 Processor, utilizing
8 GB of RAM memory. The semidefinite program was formulated
and solved with CVX  \cite{cvx}, \cite{GrantBoydYe:2006}
using the SDPT3, \cite{sdpt3}, solver and the total computation time is less than 11 seconds.
The following figure depicts  the altitude trajectory for 1000 realizations of the uncertain factors, where the system evolves in accordance to the performance specifications.
\begin{figure}[!h]
	\centering
	\includegraphics[width=140mm]{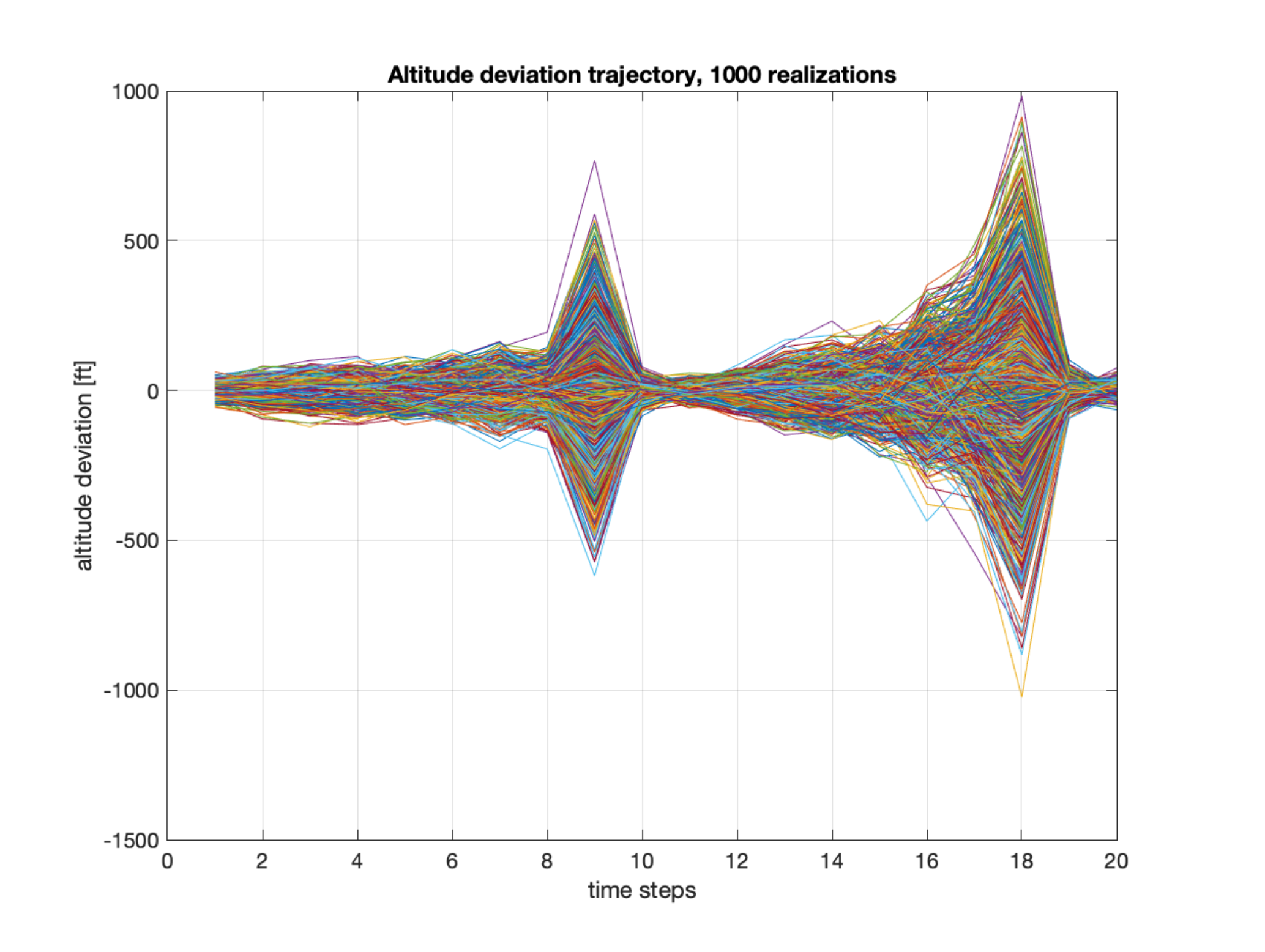}
\end{figure}

\newpage

\section{Appendix}

\subsection{Justifying claims in Section \ref{processingspec}}
\subsubsection{Justifying claim A}\label{verifyA}
In the notation of item {\bf A} and setting $\cZ=\mZ$, we should prove that $\chi$ satisfies (\ref{constraintA}), or, which is the same by (\ref{have1}), the constraint
\begin{equation}\label{cons1}
\forall (\bfzeta\in \cZ): [\bfzeta;1]^\T  \bfZ^\T [\chi]\mathcal{A}\bfZ[\chi][\bfzeta;1]-2\beta^\T\cA\bfZ[\chi][\bfzeta;1]+\Tr(\Pi^{1/2}\bfE^\T [\chi]\mathcal{A}\bfE[\chi]\Pi^{1/2})+r\leq\gamma,
\end{equation}
if and only if $\chi$ can be augmented to the feasible solution to (\ref{systemA}). \\
1) Assume that $\chi$ satisfies (\ref{cons1}), and let us define $X,x,\alpha,\delta$ according to
$$
\begin{array}{ll}
\left[\begin{array}{c|c}X&x\cr\hline x^\T &\alpha\cr\end{array}\right]=\bfZ^\T [\chi]\mathcal{A}\bfZ[\chi]&(a)\\
\delta=\max\limits_{z\in\cZ} \left[z^\T Xz+2x^\T z-2\beta^\T \cA\bfZ[\chi][z;1]\right]&(b)\\
\end{array}
$$
$(a)$ ensures the validity of (\ref{systemA}.$a$), and $(a,b)$ together with  (\ref{cons1}) ensure the validity (\ref{systemA}.$b$). Finally, identity (\ref{ident1}) and $(b)$ say that
$$
\delta=\max_{z\in\cZ}\left[z^\T Xz+2[x+p[\chi]]^\T  z+q[\chi]\right],
$$
implying the validity of (\ref{systemA}.$c$)
\\
2) Now assume that $\chi,X,x,\alpha,\delta$ satisfy (\ref{systemA}). Then
$$
\delta\geq \max_{z\in\cZ} \left[z^\T Xz+2x^\T z-2\beta^\T \cA\bfZ[\chi][z;1]\right]
$$ by (\ref{systemA}.$c$) combined with (\ref{ident1}), whence
$$
\max_{z\in\cZ}\left[[z;1]^\T\left[\begin{array}{c|c}X&x\cr\hline x^\T&\alpha\cr\end{array}\right][z;1]-2\beta^\T \cA\bfZ[\chi][z;1]\right]\leq\delta+\alpha,
$$
which combines with (\ref{systemA}.$a$) to imply that
$$
\max_{z\in\cZ}\left[[z;1]^\T \bfZ^\T[\chi]\cA\bfZ[\chi][z;1]-2\beta^\T \cA\bfZ[\chi][z;1]\right]\leq\delta+\alpha.
$$
The resulting inequality combines with (\ref{systemA}.$b$) to imply the validity of (\ref{cons1}). $\Box$
\subsubsection{Justifying claim {\bf B}}\label{verifyAA}
In view of (\ref{ident2}) and (\ref{ellitop}), constraint (\ref{constraintAA}) reads
$$
\gamma-q[\chi]\geq \Opt[\chi]:=\max_{z\in\mZ} 2z^\T p[\chi]=
\max_{y,t}\left\{ 2y^\T P^\T p[\chi]:y^\T S_ky\leq \rho t_k,k\leq K,t\in\cT\right\}.
$$
Optimization problem specifying $\Opt[\chi]$ clearly is convex, bounded and satisfies Slater condition. Applying Lagrange Duality Theorem, we get the first equivalence in the following chain:
$$
\begin{array}{ll}
&\gamma-q[\chi]\geq \Opt[\chi]\\
\Leftrightarrow&\exists \lambda\geq 0: \gamma-q[\chi]\geq \underbrace{\sup_{t\in\cT,y}\left\{ 2y^\T P^\T p[\chi]-\sum_k \lambda_k(y^\T S_ky-\rho t_k)\right\}}_{=\rho\phi_{\cT}(\lambda)+\sup_y\left\{2y^\T P^\T p[\chi]-y^\T [\sum_k\lambda_kS_k]y\right\}}\\
\Leftrightarrow&\exists \lambda\geq0,\alpha: \gamma-q[\chi]\geq \rho\phi_{\cT}+\alpha\ \&\ \left[\begin{array}{c|c}\alpha&-p^\T[\chi]P\cr\hline -P^\T p[\chi]&\sum_k\lambda_kS_k\cr\end{array}\right]\succeq0,\\
\end{array}
$$
and {\bf B} follows. $\Box$
\subsubsection{Justifying claim $\mathbf{B'}$}\label{verifyBprime}
\paragraph{Preliminaries.} We start with the well-known ``conic extension'' of Lagrange Duality Theorem:
\begin{theorem}\label{theapp} Let  $\bK$ be regular (i.e., closed convex pointed and with a nonempty interior) cone in some Euclidean space $E$, $\bK^*$ be the cone dual to $\bK$, let $W$ be a convex set in some $\bbR^n$, let $f(w):W\to\bbR$ be convex,
and let $g(w):W\to E_k$ be $\bK$-convex, meaning that $\lambda g(u)+(1-\lambda)g(v)-g(\lambda u+(1-\lambda)v)\in\bK$ for all $u,v\in W$ and $\lambda\in[0,1]$.
Consider the feasible convex optimization problem
$$
\Opt(P)=\min_w\left\{f(w): g(w)\leq_\bK 0, w\in W\right\},\eqno{(P)}
$$
where $a\geq_\bK b$ ($\Leftrightarrow b\leq_\bK a$) means that $a-b\in\bK$. Along with $(P)$, consider its ``conic'' Lagrange function
$$
L(w,\Lambda)=f(w)+\langle \Lambda,g(w)\rangle: W\times \bK^*\to\bbR
$$
and  ``conic'' Lagrange dual
$$
\Opt(D)=\max_{\Lambda}\left\{f_*(\Lambda):=\inf_{w\in W}L(w,\lambda),\Lambda\in\bK^*\right\}.\eqno{(D)}
$$
Then $\Opt(D)\leq\Opt(P)$. Moreover, assuming
that $(P)$ is bounded  and satisfies Slater condition:
$$
-g(\bar{w})\in\inter \bK
$$
for some $\bar{w}\in W$, $(D)$ is solvable and $\Opt(P)=\Opt(D)$.
\end{theorem}
Note that the standard Lagrange Duality Theorem is the special case of Theorem \ref{theapp} corresponding to the case when $\bK$ is the nonnegative orthant in some $\bbR^n$.
\par
To make the paper self-contained, here is the proof of Theorem \ref{theapp}. First, it is immediately seen that when $\Lambda\in\bK^*$, $L(w,\Lambda)$ underestimates $f(w)$ on  the feasible set of $(P)$, whence $f_*(\Lambda)\leq\Opt(P)$ on the feasible domain of $(D)$, implying that
$\Opt(D)\leq\Opt(P)$. Now assume that $(P)$ is bounded and satisfies Slater condition, and  let
$$
S=\{(\tau,x)\in\bbR\times E:\tau<\Opt(P),x\leq_\bK 0\},\, T=\{(\tau,x)\in\bbR\times E: \exists w\in W: f(w)\leq \tau, g(w)\leq_{\bK} x\}.
$$
It is immediately seen that $S$ and $T$ are nonempty nonintersecting convex sets, implying that there exists a nonzero $(\alpha,M)\in\bbR\times E$ such that
\begin{equation}\label{separ}
\sup_{(\tau,x)\in S}\left[\alpha \tau+\langle x,M\rangle\right]\leq \inf_{(\tau,x)\in T}\left[\alpha \tau+\langle x,M\rangle\right].
\end{equation}
By this relation, the right hand side quantity is finite, and since $\bbR_+\times\bK$ clearly is contained in the recessive cone of $T$, we conclude that $\alpha>0$ and $M\in\bK^*$.
As a result,  by construction of $S$, $T$
(\ref{separ}) reads
\begin{equation}\label{separ1}
\alpha\Opt(P)\leq \inf_{w\in W}\left[\alpha f(w)+\langle M,g(w)\rangle\right].
\end{equation}
When $\alpha>0$, setting $\Lambda=M/\alpha$, we get $\Lambda\in\bK^*$, and (\ref{separ1}) says that $f_*(\Lambda)\geq \Opt(P)$; since, as we have seen, $\Opt(D)\leq\Opt(P)$, we conclude that $\Lambda$ is an optimal solution to $(D)$ and $\Opt(P)=\Opt(D)$.\par It remains to verify that we indeed have $\alpha>0$. Assuming this is not the case, $\alpha=0$, and (\ref{separ1}) combines with $\bar{w}\in W$ to imply that $\langle g(\bar{w}),M\rangle\geq0$. Since $-g(\bar{w})\in\inter\bK$ and $M\in\bK^*$, we conclude that $M=0$, whence $(\alpha,M)=0$, which is not the case. $\Box$
\paragraph{Theorem \ref{theapp}$\Rightarrow$ $\mathbf{B'}$:}
In view of (\ref{ident2}) and (\ref{spectratop}), constraint (\ref{constraintAAA}) reads
$$
\gamma-q[\chi]\geq \max_{z\in\mZ} 2z^\T p[\chi]=-\Opt[\chi],\,\Opt[\chi]=
\min_{y,t}\left\{ -2y^\T P^\T p[\chi]:S_k^2[y]- \rho t_kI_{f_k}\preceq0,k\leq K,t\in\cT\right\}.
$$
Setting $E=\bbS^{f_1}\times...\times\bbS^{f_K}$, $\bK=\bbS^{f_1}_+\times...\times\bbS^{f_K}_+$, $g(y)=\Diag\{S_1^2[y]- \rho t_1I_{f_1},...,S_K^2[y]- \rho t_KI_{f_K} \}$, (\ref{constraintAAA}) becomes
\begin{equation}\label{becomes}
\begin{array}{c}
\gamma-q[\chi]\geq -\Opt[\chi],\,\Opt[\chi]=\min_{w\in W}\left\{f(y):= -2y^\T P^\T p[\chi]:g(y)\leq_\bK0\right\},\\
W=\{(y,t):t\in\cT\}.\\
\end{array}
\end{equation}
The problem specifying $\Opt[\chi]$ is of the form $(P)$ considered in Theorem \ref{theapp} and clearly is below bounded and satisfies Slater condition. Applying Theorem \ref{theapp} and noting that, as is immediately seen, $\Tr(S_k^2[y]\Lambda_k)=y^T\cS_k^*[\Lambda_k]y$, we conclude that
{\small $$
\begin{array}{l}
-\Opt[\chi]\\
\quad=-\max\limits_{\Lambda}\left\{f_*(\Lambda):=\inf\limits_{t\in\cT,y}\left[-2y^\T P^\T p[\chi]+y^T\left[\sum_k\cS_k^*[\Lambda_k]\right]y-\rho\sum_k\Tr[\Lambda_k]t_k: \Lambda=\{\Lambda_k\in\bbS^{f_k}_+,k\leq K\}\right]\right\}\\
\quad=-\max\limits_{\Lambda}\left\{\inf\limits_y\left[-2y^\T P^\T p[\chi]+y^T\left[\sum_k\cS_k^*[\Lambda_k]\right]y\right]-\rho\phi_{\cT}(\lambda[\Lambda]): \Lambda=\{\Lambda_k\in\bbS^{f_k}_+,k\leq K\}\right\}\\
\quad=\min\limits_{\Lambda}\left\{\rho\phi_{\cT}(\lambda[\Lambda])+\sup\limits_y\left[2y^\T P^\T p[\chi]-y^T\left[\sum_k\cS_k^*[\Lambda_k]\right]y\right]:\Lambda=\{\Lambda_k\in\bbS^{f_k}_+,k\leq K\}\right\}\\
\quad=\min\limits_{\Lambda,\alpha}\left\{\rho\phi_{\cT}(\lambda[\Lambda])+\alpha:\left[\begin{array}{c|c}\alpha&-p^\T[\chi]P\cr\hline
-P^\T p[\chi]&\sum_k\cS_k^*[\Lambda_k]\cr\end{array}\right]\succeq0, \Lambda=\{\Lambda_k\in\bbS^{f_k}_+,k\leq K\}\right\}\\
\end{array}
$$}\noindent
and the maxima, and therefore the minima, in the right hand side  are achieved. This conclusion combines with (\ref{becomes}) to imply $\mathbf{B'}$. $\Box$

\subsection{Proofs of Theorems \ref{thean1}, \ref{thean2}}
Theorems \ref{thean1}, \ref{thean2} are essential extensions of the results of \cite{roos} where, in retrospect, a special ellitope -- intersection of concentric ellipsoids -- was considered;

\subsubsection{Proof of Theorem \ref{thean2}}\label{appextSlemma}
\paragraph{0$^0$.} By replacing $\cT$ with $\rho\cT$, we can reduce the situation to the one where $\rho=1$, which is assumed from now on.
\paragraph{1$^0$}
For $s>0$, let
{\small $$
\begin{array}{c}
	\mathcal{X}[s]=\{x\in\bbR^n: \exists z \in \mathcal{Z}[s]:=\left\{z\in\bbR^N:\exists t\in\mathcal{T}: S_k^2[z]\leq s t_kI_{f_k},k\leq K\right\}: x=Pz\},\\
S_k[z]=\sum_iz_iS^{ki},\,S^{ki}\in\bbS^{f_k}\\
\end{array}
	$$}\noindent
be a single-parametric family of spectratopes.
Let also
$$
F(x)=x^\T Ax+2b^\T x,\,\,\Opt_*[s]=\max_x\left\{F(x): x\in \mathcal{X}[s]\right\}.
$$
Note that
$$
\Opt_*[s]=\max_z\left\{G(z)=z^\T Qz+2q^\T z: z\in\mathcal{Z}[s]\right\},\,\,Q=P^\T AP,\,q=P^\T b.
$$
Let us set
$$
Q_+= \left[\begin{array}{c|c}Q & q \cr\hline q^\T  &\cr\end{array}\right].
$$
For $\Lambda=\{\Lambda_k\in\bbS^{f_k},k\leq K\}$, let
$$
		\begin{array}{l}
		\lambda[\Lambda] = [\Tr(\Lambda_1);...;\Tr(\Lambda_K)],\\
		\mathcal{S}_k[Q]=\sum_{i,j}Q_{ij}S^{ki}S^{kj}=\sum_{i,j}Q_{ij}{1\over 2}[S^{ki}S^{kj}+S^{kj}S^{ki}]: \bbS^N\to\bbS^{f_k}\\
		\mathcal{S}_k^*[\Lambda_k]=\left[\Tr(\Lambda_k S^{ki}S^{kj})\right]_{i.j}=\left[\Tr(\Lambda_k{1\over 2}[S^{ki}S^{kj}+S^{kj}S^{ki}])\right]_{i,j}:\bbS^{f_k}\to\bbS^N.\\
		\end{array}
		$$
Note that the mappings $\mathcal{S}_k$, $\mathcal{S}_k^*$ are conjugate to each other:
		$$
		\Tr(\mathcal{S}_k[Q]\Lambda_k)=\Tr(Q\mathcal{S}_k^*[\Lambda_k]),\,\,Q\in\bbS^N,\Lambda_k\in\bbS^{f_k}
		$$
\paragraph{2$^0$.} Our local goal is to prove the following
\begin{proposition}\label{extSlemma}
	Let
	\begin{equation}\label{eqS200}
	\Opt[s]=\min_{\Lambda,\mu}\left\{s\phi_{\mathcal{T}}(\lambda[\Lambda])+\mu:\begin{array}{l}
	\Lambda=\{\Lambda_k\succeq0,k\leq K\},\mu\geq 0\\
	Q_+\preceq \left[\begin{array}{c|c}\sum_k\mathcal{S}_k^*[\Lambda_k]&\cr\hline &\mu\cr\end{array}\right]\\
	\end{array}\right\}.
	\end{equation}
	Then
	\begin{equation}\label{eqS201}
	\Opt_*[s]\leq \Opt[s]\leq \Opt_*[\varkappa s]
	\end{equation}
	with
	\begin{equation}\label{eqS202}
	\varkappa=2\ln\left(8\sum_{k=1}^Kf_k\right).
	\end{equation}
\end{proposition}
\noindent
\textbf{Proof of Proposition.} Replacing $\mathcal{T}$ with $s\mathcal{T}$, we can assume once for ever that $s=1$. Now, setting
$$
\mathbf{T}=\cl\{[t;\tau]: \tau>0, t/\tau\in\mathcal{T}\}
$$
we get a closed convex pointed cone with a nonempty interior, the dual cone being
$$
\mathbf{T}_*=\{[g;s]: s\geq\phi_{\mathcal{T}}(-g)\}.
$$
Problem (\ref{eqS200}) with $s=1$ is the conic problem
$$
\Opt[1]=\min_{\Lambda,\tau,\mu}\left\{\tau+\mu:\begin{array}{l}\Lambda_k\succeq0,k\leq K,\mu\geq 0\\
{[-\lambda[\Lambda];\tau]}\bar{}\in\mathbf{T}_*\\
\left[\begin{array}{c|c}\sum_k\mathcal{S}_k^*[\Lambda_k]&\cr\hline &\mu\cr\end{array}\right]-Q_+\succeq0\\
\end{array}\right\}.\eqno{(*)}
$$
It is immediately seen that this problem is strictly feasible and solvable, so that $\Opt[1]$ is the optimal value in the  conic dual of $(*)$. To get the latter problem, let $P_k\succeq0$, $p\geq0$, $[t;s]\in\mathbf{T}$, and
$W=\AR{cc} {V&v \\ v^\T &w} \succeq0$
 be Lagrange multipliers for the respective constraints in $(*)$. Aggregating the constraints with these weights, we get
$$
\sum_k\Tr(P_k\Lambda_k)+p\mu-t^\T \lambda[\Lambda]+\tau s +\Tr(V\sum_k\mathcal{S}_k^*[\Lambda_k])-\Tr(VQ)-2v^\T q+w\mu\geq 0,
$$
that is,
$$
\sum_k\Tr([P_k+\mathcal{S}_k[V]-t_kI_{d_k}]\Lambda_K) +[p+w]\mu +\tau s \geq \Tr(VQ)+2v^\T q
.$$
To get the dual problem, we add to the above restrictions on the Lagrange multipliers the requirement that the left hand side in the latter inequality, identically in $\Lambda_k$, $\tau$, and $\mu$ should be equal to the objective in $(*)$, implying, in particular, that $t\in\mathcal{T}$, $s=1$, and $w\leq 1$. We see that the dual problem, after evident simplifications, becomes
\begin{equation}\label{eqS500}
\Opt[1]=\max_{V,v,t}\left\{\Tr(VQ)+2v^\T q:\begin{array}{l} t\in\mathcal{T},\,\mathcal{S}_k[V]\preceq t_kI_{f_k},k\leq K,\\
\left[\begin{array}{c|c}V&v\cr\hline v^\T &1\cr\end{array}\right]\succeq0.
\\
\end{array}\right\}
\end{equation}
By definition, $$
\Opt_*[1]=\max_{z,t}\left\{z^\T Qz+2q^\T z:t\in\mathcal{T}, S_k^2[z]\preceq t_kI_{f_k},k\leq K\right\}.
$$
If $(z,t)$ is a feasible solution to the latter problem, then $V=zz^\T ,v=z,t$ is a feasible solution to (\ref{eqS500})  with the same value of the objective, implying that
$\Opt[1]\geq\Opt_*[1]$, as stated in the first relation in (\ref{eqS201}) (recall that we are in the case of $s=1$). \par
We have already stated that problem (\ref{eqS500}) is solvable. Let $V_*,v_*,t^*$ be its optimal solution, and let
$$
X=\left[\begin{array}{c|c}V_*&v_*\cr\hline v_*^\T &1\cr\end{array}\right].
$$
Let $\chi$ be Rademacher random vector\footnotemark \footnotetext{random vector with independent entries taking values $\pm1$ with probabilities 1/2.} of the same size as the one of $Q_+$, let
$$
X^{1/2}Q_+X^{1/2}=U\Diag\{\lambda\}U^\T
$$
with orthogonal $U$, and let
$\zeta=X^{1/2}U\chi=[\xi;\tau]$, where $\tau $ is the last entry in $\zeta$. Then
\begin{equation}\label{eqS400}
\begin{array}{l}
\zeta^\T Q_+\zeta=\chi^\T U^\T X^{1/2}Q_+X^{1/2}U\chi=\chi^\T \Diag\{\lambda\}\chi\\
=\sum_i\lambda_i=\Tr(X^{1/2}Q_+X^{1/2})=\Tr(XQ_+)=\Opt[1].\\
\end{array}
\end{equation}
Next, we have
$$
\zeta\zeta^\T =\left[\begin{array}{c|c}\xi\xi^\T &\tau\xi\cr\hline \tau\xi^\T &\tau^2\cr\end{array}\right]=X^{1/2}U\chi\chi^\T U^\T X^{1/2},
$$
whence
\begin{equation}\label{eqS220}
\bbE\{\zeta\zeta^\T \}=\left[\begin{array}{c|c}\bbE\{\xi\xi^\T \}&\bbE\{\tau\xi\}\cr\hline \bbE\{\tau\xi^\T \}&\bbE\{\tau^2\}\\
\end{array}\right]=X=\left[\begin{array}{c|c}V_*&v_*\cr\hline v_*^\T &1\cr\end{array}\right].
\end{equation}
In particular,
$$
\mathcal{S}_k[\bbE\{\xi\xi^\T \}]=\mathcal{S}_k[V_*]\preceq t_k^*I_{f_k},\,\,k\leq K.
$$
We have $\xi=W\chi$ for certain rectangular matrix $W$ such that $V_*=\bbE\{\xi\xi^\T \}=
\bbE\{\ W \chi\chi^\T  W^\T \}=W W^\T $. Consequently,
$$
S_k[\xi]=S_k[W\chi]=\sum_i\chi_i\bar{S}^{ki}:=\bar{S}_k[\chi]
$$
with some symmetric matrices $\bar{S}^{ki}$. As a result,
$$
\mathcal{S}_k[\xi\xi^\T ]=S_k^2[\xi]=\sum_{i,j}\chi_i\chi_j\bar{S}^{ki}\bar{S}^{kj},
$$
and taking expectation we get
$$
\sum_i[\bar{S}^{ki}]^2=\mathcal{S}_k[\bbE\{\xi\xi^\T \}]=\mathcal{S}_k[V_*]\preceq t_k^*I_{f_k},
$$
whence by Noncommutative Khintchine Inequality (see, e.g., \cite[Theorem 4.45]{A_A_20}) one has
\begin{equation}
\label{eqS210}
\Pr\{S_k^2[\xi]\preceq t^*_k r^{-1}I_{f_k}\}=\Pr\{\bar{S}_k^2[\chi]\preceq t^*_kr^{-1}I_{f_k}\}\leq 1 - 2f_k\exp\{-{1\over 2r}\},\,0<r\leq 1.
\end{equation}
Next, invoking (\ref{eqS220}) we have
\begin{equation}
\label{eqS221}
\bbE\{\tau^2\}=\bbE\{[\zeta\zeta^\T ]_{N+1,N+1}\}=1
\end{equation}
and
$$
\tau=\beta^\T \chi
$$
for some vector $\beta$ with $\|\beta\|_2=1$ due to (\ref{eqS221}). Now let us invoke the following fact \cite[Lemma A.1]{roos}
\begin{lemma}
	\label{lemSlem}
	Let $\beta$ be a deterministic $\|\cdot\|_2$-unit vector in $
	\bbR^N$ and $\chi$ be $N$-dimensional Rademacher random vector. Then $\Pr\{|\beta^\T \chi|\leq 1\}\geq 1/3$.
\end{lemma}

Now let $\bar{r}$ be given by the relation
$$
\bar{r}={1\over 2\ln(8\sum_kf_k)},
$$
implying that
$$
2\sum_{k=1}^K f_k\exp\{-{1\over 2\bar{r}}\}<1/3.
$$
By (\ref{eqS210}) and Lemma \ref{lemSlem} there exists a realization $\bar{\zeta}=[\bar{\xi};\bar{\tau}]$ of $\zeta$ such that
\begin{equation}\label{eqS300}
S_k^2[\bar{\xi}]\preceq t^*_k \bar{r}^{-1}I_{f_k},\,k\leq K, \ \&\ |\bar{\tau}|\leq1.
\end{equation}
Invoking (\ref{eqS400}) and taking into account that $|\bar{\tau}|\leq1$ we have
$$
\Opt[1]=\bar{\zeta}^\T Q_+\bar{\zeta}=\bar{\xi}^\T Q\bar{\xi}+2\bar{\tau}q^\T \bar{\xi}\leq \widehat{\xi}^\T Q\widehat{\xi}+2q^\T \widehat{\xi},
$$
where $\widehat{\xi}=\bar{\xi}$ when $q^\T \bar{\xi}\geq0$ and $\widehat{\xi}=-\bar{\xi}$ otherwise. In both cases from the first relation in (\ref{eqS300})
we conclude that $\widehat{\xi}\in\mathcal{Z}[1/\bar{r}]$, and we arrive at $\Opt[1]\leq \Opt_*[2\ln(8\sum_kf_k)]$. Proposition \ref{extSlemma} is proved.
\paragraph{3$^0$.} We are ready to complete the proof of Theorem\ref{thean2}. Indeed, by replacing $\cT$ with $\rho\cT$, we reduce the situation considered in Theorem to the one where $\rho=1$. Now,
using the notation from Theorem,  let us set
$$
n=n_\zeta,\,N=\bar{n}, Q_+=\left[\begin{array}{c|c}X&x\cr\hline x^\T &\cr\end{array}\right].
$$
thus ensuring that $\cX[1]=\mathscr{D}[1]$. $Q_*$, $K$ and the entities participating in Theorem \ref{thean2} give rise to the data participating in Proposition \ref{extSlemma}. In terms of this Proposition, the fact that $X,x,\xi$ is feasible for $(C[s])$ with $s>0$ is nothing but the relation
$$
\Opt_*[s]\leq\xi,
$$
and the fact that $X,x,\xi$ can be extended by properly selected $\Lambda,\tau$ to a feasible solution to $(\bS[s])$ is exactly the relation
$$
\Opt[s]\leq\xi.
$$
As a result, Theorem \ref{thean2} becomes an immediate consequence of Proposition \ref{extSlemma}, with item $(i)$ of Theorem readily given by the left, and item $(ii)$ -- by the right inequality in (\ref{eqS201}).

\subsubsection{Proof of Theorem \ref{thean1}}\label{appextSlemmaEll}
\paragraph{1$^0$.} Consider a single-parametric family of ellitopes
\begin{equation}\label{postellit}
	\begin{array}{c}
	\cX[s]=\{x:\exists z\in\cZ[s]: x=Pz\},\,\mathcal{Z}[s]=\{z\in\bbR^N:\exists t\in\mathcal{T}: z^\T S_kz\leq s t_k,k\leq K\}\\
	\end{array}, s>0.
	\end{equation}
Let also
$$
F(x)=x^\T Ax+2b^\T x,\,\,\Opt_*[s]=\max_x\left\{F(x): x\in \mathcal{X}[s]\right\}/
$$
Note that
$$
\Opt_*[s]=\max_z\left\{G(z)=z^\T Qz+2q^\T z: z\in\mathcal{Z}[s]\right\},\,\,Q=P^\T AP,\,q=P^\T b.
$$
Let us set
$$
Q_+= \left[\begin{array}{c|c}Q & q \cr\hline q^\T  &\cr\end{array}\right].
$$
The ellitopic version of Proposition \ref{extSlemma}
reads as follows:
\begin{proposition}\label{extSlemmaEll}
	Let
	\begin{equation}\label{eqS200Ell}
	\Opt[s]=\min_{\lambda,\mu}\left\{s\phi_{\mathcal{T}}(\lambda)+\mu:\begin{array}{l}
	\lambda=\{\lambda_k\geq0,k\leq K\},\mu\geq 0\\
	Q_+\preceq \left[\begin{array}{c|c}\sum_k\mathcal{S}_k^*[\Lambda_k]&\cr\hline &\mu\cr\end{array}\right]\\
	\end{array}\right\}.
	\end{equation}
	Then
	\begin{equation}\label{eqS201Ell}
	\Opt_*[s]\leq \Opt[s]\leq \Opt_*[\varkappa s]
	\end{equation}
	with
	\begin{equation}\label{eqS202Ell}
	\varkappa=3\ln(6K).
	\end{equation}
\end{proposition}
Theorem \ref{thean1} can be immediately derived from this Proposition in exactly the same fashion as Theorem \ref{thean2} above was derived from Proposition \ref{extSlemma}.
\paragraph{3$^0$. Proof of Proposition \ref{extSlemmaEll}} is similar to the one of Proposition \ref{extSlemma}. Specifically,
\paragraph{3$^0$.1.} We have
		$$
		\Opt_*[s]=\max_{z,t}\left\{z^\T Qz+2q^\T z:t\in\mathcal{T}, z^\T S_kz\leq s t_k,k\leq K\right\},
		$$
		and
		$$
		\begin{array}{rcl}
		\Opt[s]&=&\min_{\lambda,\mu}\left\{s\phi_\mathcal{T}(\lambda)+\mu:\begin{array}{l}\lambda\geq0,\mu\geq0\\
		\left[\begin{array}{c|c}\sum_k\lambda_kS_k&\cr\hline &\mu\cr\end{array}\right]-Q_+\succeq 0\\
		\end{array}\right\}\\
		&=&\min_{\lambda,\mu,\tau}\left\{s\tau+\mu:\begin{array}{l}\lambda\geq0,\mu\geq0,[-\lambda;\tau]\in\mathbf{T}_*\\
		\left[\begin{array}{c|c}\sum_k\lambda_kS_k&\cr\hline &\mu\cr\end{array}\right]-Q_+\succeq 0\\
		\end{array}\right\}\\
		\end{array}.
		$$
\paragraph{3$^0$.2.} Same as in the spectratopic case, it suffices to consider the case when $s=1$, which is assumed from now on.\par
The problem dual to the conic problem specifying $\Opt[1]$ reads
		\begin{equation}\label{eqS610Ell}
		\Opt[1]=\max_{V,v,t}\left\{\Tr(QV)+2q^\T v:\begin{array}{l}t\in\mathcal{T},\Tr(S_kV)\leq t_k,k\leq K\\
		\left[\begin{array}{c|c}V&v\cr\hline v^\T &1\cr\end{array}\right]\succeq0\\
		\end{array}\right\}
		\end{equation}
Every feasible solution $(z,t)$ to the problem specifying $\Opt_*[1]$ gives rise to the feasible solution $(V=zz^\T ,v=z,t)$ to (\ref{eqS610Ell}) with the same value of the objective,
implying that $\Opt_*[1]\leq\Opt[1]$
\paragraph{3$^0$.3.} It is easily seen that (\ref{eqS610Ell}) has an optimal solution $V_*,v_*,t^*$.
		Setting
		$$ X=\left[\begin{array}{c|c}V_*&v_*\cr\hline v_*^\T &1\cr\end{array}\right],$$
		we set
		$$
		X^{1/2}Q_+X^{1/2}=U\diag\{\lambda\}U^\T
		$$
		with orthogonal $U$, set
		$$
		\zeta=X^{1/2}U\chi
		$$
		with Rademacher random $\chi$, thus ensuring that
		{\small\begin{equation}\label{eqS660Ell}
			\zeta^\T Q_+\zeta=\chi^\T U^\T X^{1/2}Q_+X^{1/2}U\chi=\chi^\T \diag\{\lambda\}\chi=\sum_i\lambda_i=\Tr(XQ_+)=\Opt[1].
			\end{equation}}\noindent
		At the same time, setting $\zeta=[\xi;\tau]$ with scalar $\tau$, we get
		{\small $$
			\left[\begin{array}{c|c}\bbE\{\xi\xi^\T \}&\bbE\{\tau\xi\}\cr\hline
			\bbE\{\tau\xi^\T \}&\bbE\{\tau^2\}\cr\end{array}\right] =\bbE\{\zeta\zeta^\T \}=\bbE\{X^{1/2}U\chi\chi^\T U^\T X^{1/2}\}=X=\left[\begin{array}{c|c}
			V_*&v_*\cr\hline v_*^\T &1\cr\end{array}\right].
			$$}\noindent
		As before, we have $\xi=W\chi$ for some deterministic matrix $W$, and $\bbE\{\xi^\T S_k\xi\}=\Tr(S_kV_*)\leq t_k^*$, whence
		$\Tr(W^\T S_kW)=\bbE\{\chi^\T W^\T S_kW\chi\}=\bbE\{\xi^\T S_k\xi\}\leq t_k^*$. By \cite[Lemma 4.28]{A_A_20} the relations $W^\T S_kW\succeq0$, $\Tr(W^\T S_kW)\leq t_k^*$  imply that
		$$
		\Pr\{\xi^\T S_k\xi>\gamma t_k^*\}=\Pr\{\chi^\T W^\T S_kW\chi>\gamma t_k^*\}\leq \sqrt{3}\exp\{-{\gamma\over 3}\},\, \forall \gamma > 0,
		$$
(cf. the proof of \cite[Proposition 4.6]{A_A_20}).
\paragraph{3$^0$.4.} Now the reasoning completely similar to the one in the spectratopic case shows that when $\varkappa=3\ln(6K)$, there exists a realization $\bar{\zeta}=[\bar{\xi};\bar{\tau}]$ of $\zeta$ such that $S_k^2[\bar{\xi}]\leq t_k^*\varkappa$ for all $k$ and $|\bar{\tau}|\leq 1$. This conclusion, taken together with (\ref{eqS660Ell}), results in $\Opt[1]\leq\Opt_*[\varkappa]$, thus completing the proof.

\section*{Acknowledgment}

The authors gratefully acknowledge the support of National Science Foundation grants 1637473, 1637474, and the NIFA grant 2020-67021-31526.

\end{document}